\documentclass[10pt]{amsart}
\usepackage{amsfonts}
\usepackage{color}
\usepackage{setspace}
\allowdisplaybreaks[4]

\definecolor{c20}{rgb}{0.,0.7,0.}
\definecolor{c30}{rgb}{0.,0.,1.}
\definecolor{c40}{rgb}{1,0.1,0.7}
\definecolor{c50}{rgb}{1,0,0}
\definecolor{c60}{rgb}{1,0.9,0.1}

\def\jE#1{\textcolor{c20}{#1}}
\def\jE#1{#1}
\def\iE#1{\textcolor{c20}{#1}}
\def\iE#1{#1}
\def\AH#1{\textcolor{c20}{#1}}

\def\bE#1{\textcolor{c20}{#1}}
\def\iH#1{\textcolor{c20}{#1}}
\def\iH#1{#1}
\def\bE#1{#1}
\def\aE#1{\textcolor{c20}{#1}}
\def\aE#1{#1}
\def\LE#1{\textcolor{c20}{#1}}
\def\LE#1{#1}

\def\TT#1{\textcolor{c50}{#1}}
\def\TT#1{#1}
\def\Ta#1{\textcolor{c50}{#1}}
\def\Ta#1{#1}
\def\Tan#1{\textcolor{c50}{#1}}
\def\Tan#1{#1}

\def\Quan#1{\textcolor{c50}{#1}}
\def\Quan#1{#1}

\def\x{\vk{x}}
\def\y{\vk{y}}

\def\M{\vk{M}}
\def\RDI{\mathfrak{R}(\delta_i)}
\def\xyrp{\x,\y_1,\y_2 \in \R^p}
\newcommand{\EE}[1]{\mathbb{E}\left\{#1 \right\}}

\newcommand{\kb}[1]{\boldsymbol{#1}}
\newcommand{\vk}[1]{\kb{#1}}

\newcommand{\ve}{\varepsilon}

\newcommand{\abs}[1]{\left\lvert #1 \right\rvert}
\newcommand{\Abs}[1]{ \biggl \lvert #1 \biggr \rvert}
\newcommand{\ABs}[1]{ \biggl \lvert #1 \biggr \rvert}

\newcommand{\E}[1]{\mathbb{E}\left\{#1\right\}}

\newcommand{\pk}[1]{\mathbb{P} \left\{ #1 \right \} }

\newcommand{\R}{\mathbb{R}}

\newcommand{\N}{\mathbb{N}}
\newcommand{\inr}{\in \R}
\newcommand{\inn}{\in \N}
\newcommand{\ldot}{,\ldots,}

\newcommand{\limit}[1]{\lim_{#1 \to   \infty}}

\newcommand{\BQN}{\begin{eqnarray}}
\newcommand{\EQN}{\end{eqnarray}}
\newcommand{\BQNY}{\begin{eqnarray*}}
\newcommand{\EQNY}{\end{eqnarray*}}

\newcommand{\BS}{\begin{sat}}
\newcommand{\ES}{\end{sat}}
\newcommand{\BT}{\begin{theo}}
\newcommand{\ET}{\end{theo}}
\newcommand{\BL}{\begin{lem}}
\newcommand{\EL}{\end{lem}}
\newcommand{\BK}{\begin{korr}}
\newcommand{\EK}{\end{korr}}

\newcommand{\BD}{\begin{de}}
\newcommand{\ED}{\end{de}}

\newcommand{\BIT}{\begin{itemize}}
\newcommand{\EIT}{\end{itemize}}
\newcommand{\BDI}{\begin{description}}
\newcommand{\EDI}{\end{description}}

\newcommand{\BRM}{\begin{remarks}}
\newcommand{\ERM}{\end{remarks}}

\newcommand{\BEL}{\begin{lem}}
\newcommand{\EEL}{\end{lem}}

\def\SL{{\mathcal S}_l}

\def\SJ{{\mathcal S}_j}

\def\RL{{\mathcal R}_l}
\def\RJ{{\mathcal R}_j}
\def\IF{\infty}

\newtheorem{theo}{Theorem}[section]
\newtheorem{sat}[theo]{Proposition}
\newtheorem{de}[theo]{Definition}
\newtheorem{lem}[theo]{Lemma}

\newtheorem{korr}[theo]{Corollary}

\newtheorem{remarks}[theo]{Remarks}

\newcommand{\nelem}[1]{{Lemma \ref{#1}}}

\newcommand{\prooftheo}[1]{ \textsc{\bf Proof of Theorem} \ref{#1}:}

\newcommand{\prooflem}[1]{\textsc{\bf Proof of Lemma} \ref{#1}:}

\newcommand{\COM}[1]{}

\newcommand{\QED}{\hfill $\Box$}

\def\cE#1{#1}

\def\IF{\infty}

\def\X{\mathbf X}

\newcommand{\todis}{\stackrel{d}{\to}}
\newcommand{\toprob}{\stackrel{p}{\to}}

\begin{document}

\title[Piterbarg's max-discretisation theorem]{Piterbarg's max-discretisation theorem for stationary vector Gaussian processes observed on different grids}

\author{Enkelejd  Hashorva}
\address{Enkelejd Hashorva, Department of Actuarial Science, Faculty of Business and Economics (HEC Lausanne), University of Lausanne,\\
UNIL-Dorigny, 1015 Lausanne, Switzerland
}
\email{Enkelejd.Hashorva@unil.ch}

\author{Zhongquan Tan}
\address{Zhongquan Tan, College of Mathematics, Physics and Information Engineering, Jiaxing University, Jiaxing 314001, PR China, and Department of Actuarial Science, Faculty of Business and Economics (HEC Lausanne), University of Lausanne,\\
UNIL-Dorigny, 1015 Lausanne, Switzerland, Corresponding author}
\email{zhongquan.tan@unil.ch}

\bigskip
\date{\today}
\maketitle
\COM{
{\small\it 1. Department of Actuarial Science, Faculty of Business and Economics, University of Lausanne,
UNIL-Dorigny, 1015 Lausanne, Switzerland. E-mail: Enkelejd.Hashorva@unil.ch}

{\small\it 2. College of Mathematics, Physics and Information Engineering, Jiaxing University, Jiaxing 314001, PR China. Corresponding author. Tel.: +86 13957395469.
E-mail: tzq728@163.com}\\
}

{\bf Abstract:} In this paper we derive Piterbarg's max-discretisation theorem \bE{for two} different  grids considering centered stationary vector Gaussian process\iH{es}. So far in the literature results in this direction have been derived
 for the joint distribution of \iE{the} \bE{maximum} of
Gaussian processes over $[0,T]$ and over \bE{a grid} $ \mathfrak{R}(\delta_1(T))=\{k\delta_1(T): k=0,1,\cdots\}$.
In this paper we extend \iH{the recent findings} by considering additionally the \bE{maximum} over another
grid $ \mathfrak{R}(\delta_2(T))$. We derive the joint limiting distribution of \bE{maximum} of stationary Gaussian vector processes for different choices of such grids \Quan{by} letting $T\to \IF$. \bE{As a by-product we find that}
the joint limiting distribution of the \bE{maximum} over different grids, which we refer to as the Piterbarg distribution,
 is in the case of weakly dependent Gaussian processes a max-stable distribution.



\section{Introduction}

Let $\{X(t), t\geq0\}$ be a centered stationary Gaussian process with  continuous sample paths,
unit variance and correlation
function $r(\cdot)$ which satisfies for some $\alpha\in(0, 2]$
\begin{eqnarray}
\label{eq1.1}
r(t)=1-C|t|^{\alpha}+o(|t|^{\alpha})\ \ \mbox{as}\ \ t\rightarrow 0\ \ \mbox{and}\ \ r(t)<1\ \ \mbox{for}\ \ t\iH{\not=0},
\end{eqnarray}
where $C$ is some positive constant.
In various applications only realisations of $X$ on a discrete time grid are possible.
For simplicity, in this paper we shall consider uniform grids of points $ \mathfrak{R}(\delta)=\{k\delta: k=0,1,\cdots\}$ where
$\delta:=\delta(T)>0$ depends on the parameter $T>0$. In view of the findings of
Berman (see \cite{Berman64,Berman92}) \aE{ the maximum of  $X$} taken over such a discrete grid has a limiting Gumbel distribution
if
\BQN \label{gridS}
\limit{T}(2\ln  T)^{1/\alpha} \delta(T)
=D, 
\EQN
with $D=\infty$ and the Berman condition
\begin{eqnarray}
\label{eq1.2}
\limit{T} r(T)\ln  T=r
\end{eqnarray}
holds for $r=0$. Specifically, for the maximum $M(\iH{\delta, T})= \max_{i: 0 \leq i\delta \leq T} X(i\delta)$ \iE{over
\bE{$ \mathfrak{R}(\delta) \cap [0,T]$}}
we have
\BQNY
\limit{T} \sup_{x\inr} \Abs{ \pk{  a_T(M(\iH{\delta, T})- b_{\delta,T}) \le x}- e^{-e^{-x}} }=0,
\EQNY
provided that both \eqref{gridS} and \eqref{eq1.2} hold, where
\BQN \label{deltaT}
a_T=\sqrt{2 \ln T}, \quad b_{\delta,T}= a_T- \frac{\ln (   a_T  \delta \sqrt{2 \pi }) }{a_T}, \quad T>0.
\EQN
\iE{For} the maximum over $[0,T]$ defined thus \aE{as} $M(T)= \max_{ t\in [0,T] } X(t)$
it is well-known (see e.g., \cite{leadbetter1983extremes, AlbinPHD, Albin1990, Berman92, Pit96}) that (\ref{eq1.1}) and (\ref{eq1.2}) imply
\BQN
\label{eq1.A}
\limit{T} \sup_{x\inr} \Abs{ \pk{  a_T(M(T)- b_{T}) \le x}- e^{-e^{-x}} }=0,
\EQN
where
\begin{eqnarray}
\label{bT}
 b_{T}=
 a_T+ a_T^{-1} \ln((2\pi)^{-1/2}C^{1/\alpha}H_{\alpha}a_T^{-1+2/\alpha})
\end{eqnarray}
and  $H_{\alpha}\in (0,\IF)$ denotes  Pickands constant, see \cite{PicandsA, Pit72, Berman82,
leadbetter1983extremes, Albin1990, Pit96, debicki2008note, AlbinC, DikerY, DebKo2013, DEJ14, HJ14c} for more details and generalisations of $H_\alpha$. \\
The seminal contribution \cite{Pit2004} derives the joint convergence as $T\to \IF$ of $M(T)$ and $M(\delta,T)$
 showing their asymptotic independence, i.e.,
\BQNY
\limit{T} \sup_{x,y\inr} \Abs{ \pk{  a_T(M(T)- b_{T}) \le x,  a_T(M(\delta,T)- b_{\delta, T}   \le y}- e^{-e^{-x} - e^{-y}} }=0.
\EQNY
Hereafter \aE{we set} $B^*_{\alpha/2}(t):= \sqrt{2}B_{\alpha/2}(t)-\iH{\abs{t}}^{\alpha},t\ge 0$ with $B_{\alpha}$ a standard fractional Brownian motion with Hurst index $\alpha/2\in (0,1)$; \iH{recall that $\delta=\delta(T)$ is given by \Quan{\eqref{gridS}}}. Define  further for any $D>0$
$$H_{D,\alpha}=\lim_{\lambda\rightarrow\infty}\lambda^{-1}
\E{e^{ \max_{k\inn:kD\in[0,\lambda]} B_{\alpha/2}^*(kD)}}
\in (0,\IF)$$
and set \aE{(the constant  $C>0$ below relates to \eqref{eq1.1})}
\BQN \label{bDT}
b_{T}(D)=a_T+a_T^{-1} \ln ((2\pi)^{-1/2}C^{1/\alpha}H_{D,\alpha}a_T^{-1+2/\alpha}).
\EQN
\iH{For}  $\mathfrak{R}(D a_T^{-2/\alpha}), D>0$  (in this case the grid  is called \iH{Pickands} grid \iH{and $\delta=\delta(T)=D a_T^{-2/\alpha}$}),
 then  in view of \cite{Pit2004}, Theorem 2  the stated asymptotic independence does not hold since
\BQNY
\limit{T} \sup_{x,y\inr} \Abs{ \pk{  a_T(M(T)- b_{T}) \le x,  a_T(M(\delta,T)- b_{T}(D)) \le y}-
e^{-e^{-x}-e^{-y}+H_{D,\alpha}^{\ln  H_{\alpha}+x,\ln  H_{D,\alpha}+y} } }=0,
\EQNY
 where the function $H_{D,\alpha}^{x,y}$ is defined for any $x,y\inr$ as
\BQN\label{HDA}
H_{D,\alpha}^{x,y}=\lim_{\lambda\rightarrow\infty}\lambda^{-1} H_{D,\alpha}^{x,y}(\lambda)\in (0,\IF),
\EQN
with
\begin{eqnarray*}
H_{D,\alpha}^{x,y}(\lambda)=\int_{s \inr}
e^{s}\pk{\max_{t\in[0,\lambda]}B_{\alpha/2}^*(t)>s+x,
\max_{k\inn:k D\in [0,\lambda]}B_{\alpha/2}^*(kD)>s+y}\, ds.
\end{eqnarray*}
\LE{\aE{Since it \iH{follows} that for any $w\inr$}
\BQN \label{eqzIF}
\lim_{x\to - \IF} H_{D,\alpha}^{x,w}=\aE{ e^{-w} H_{D,\alpha}, \quad  \lim_{y\to - \IF} H_{D,\alpha}^{w,y} = e^{-w} H_{\alpha}} \in (0,\IF),
\EQN
then
$$ Q(x,y)=e^{-e^{-x}-e^{-y}+H_{D,\alpha}( \ln  H_{\alpha}+x,\ln  H_{D,\alpha}+y)}  , \quad x,y\inr$$
is a bivariate distribution function which  has Gumbel marginals $Q(z,\IF)=Q(\IF,z)=e^{-e^{-z}}, z\inr$.}
Moreover $Q$ is a bivariate max-stable distribution, which we shall refer to as Piterbarg distribution. This multivariate distribution
is of some independent interest for statistical modelling of dependent multivariate risks. \\
 In the extreme case of a dense grid, which in the terminology of \cite{Pit2004} \bE{means}
that \eqref{gridS} holds for $D=0$, then by Theorem 3 in \cite{Pit2004}
\BQNY
\limit{T} \sup_{x,y\inr} \Abs{ \pk{  a_T(M(T)- b_{T}) \le x,  a_T(M(\iH{\delta, T})- b_{T}) \le y}- e^{-e^{-\min(x,y)}} }=0
\EQNY
thus the continuous time and \bE{the} discrete time maxima are asymptotically completely dependent.\\
\iH{In case of} two different uniform girds $\mathfrak{R}(\delta_1)$ and $\mathfrak{R}(\delta_2)$ a natural question
\iH{that arises} is:\\
What is the
joint limiting behaviour of $M(T), M(T,\delta_1), M(T,\delta_2)$ \iH{for different types of grids}? \\
\iH{Motivated by this question, our findings this contribution include:}\\
a) \aE{We  show that $M(T,\delta_1)$ and $M(T,\delta_2)$ are always asymptotically independent if one grid is sparse and the other grid is
Pickands or  dense. Further, we obtain the joint limiting distribution if one of the grids is Pickands, and the other grid is
 Pickands or  dense.} \\
b) The Berman condition is relaxed \bE{by} assuming that \eqref{eq1.2} holds \iH{for some}  $r\in [0,\IF)$. When $r>0$ the \iH{Gaussian} process $X$ is said to be strongly dependent, see \cite{mittal1975limit, Pit96, Nadarajah09, tanPengH2012, tanH2012, DebSoja} for details on the extremes of such Gaussian processes.  The contribution \cite{TanYuebStat} derives  Piterbarg's max-discretisation theorem for strongly dependent Gaussian processes.
 In applications, often modelling of \iE{the} maximum of \iE{functionals} of a Gaussian vector process is of interest, see e.g., \cite{MR1747100,MR2654766,DebickiHJminima}. Our results in this paper are derived for the more general framework of Gaussian vector processes extending the recent findings of \cite{TanH2014} by considering simultaneously two different grids.
 This paper \AH{highlights} the role of different grids in the approximation of the maximum over a continuous interval. Our results are therefore of interest for simulation studies, which was the main motivation of \cite{Pit2004,Huesler2004,HPit2004,TurkmanA, TurkmanB, tan2012piterbarg, Nadarajah13}.
c) As a by-product \aE{we show that for weakly dependent stationary Gaussian processes the limiting distributions
are max-stable.} In Extreme Value Theory max-stable distributions and processes are characterised in different ways, see e.g., \cite{Faletal2010,DiekerMik}. In order for a \iE{multivariate} max-stable  distribution to be also useful \bE{for statistical} modelling, it is important to find how \bE{that} distribution \iE{approximates
the maxima} of certain sequences (or triangular arrays).  Piterbarg max-stable distributions
are therefore important since we show also their usefulness in the approximations of maxima \iE{over} different grids.


Organisation of the article is as \bE{follows}. Our main results are presented in the next section.
All the proofs are relegated to Section \aE{3} which is followed  \iE{by an} Appendix.

\section{Main results}
We shall investigate  in the following the asymptotics of maxima over different grids of a centered stationary multivariate $p$-dimensional
Gaussian process $\{\vk{X}(t), t\ge 0\}$. Each component $X_k,k\le p$ of $\vk{X}$ is assumed to have a constant variance function equal to 1,
continuous sample paths  and correlation function $r_{kk}(t)=Cov(X_{k}(s),X_{k}(s+t))$ which satisfies for any index $k\le p$
\begin{eqnarray}
\label{eq1.3}
r_{kk}(t)=1-C|t|^{\alpha}+o(|t|^{\alpha})\ \ \mbox{as}\ \ t\rightarrow 0\ \ \mbox{and}\ \ r_{kk}(t)<1\ \ \mbox{for}\ \ \iH{t\not=0}
\end{eqnarray}
for some positive constants $C$.
Hereafter we suppose that
$\vk{X}$ has jointly stationary components with cross-correlation function $r_{kl}(t)=Cov(X_{k}(s),X_{l}(s+t))$ which does not depend on $s$ for any $s,t$ positive. The strong dependence condition for the vector \iE{Gaussian} process $\vk{X}$ reads 
\begin{eqnarray}
\label{eq1.4}
\lim_{T\to \infty} r_{kl}(T) \ln T= r_{kl}\in [0,\infty), \quad 1\leq k,l\leq p.
\end{eqnarray}
\iE{In} order to exclude the possibility that $\abs{X_{k}(t)}= \abs{X_{l}(t+t_{0})}$ for some $k\neq l, t_0 >0$  \begin{eqnarray}
\label{eq1.5}
\max_{k\neq l}\sup_{t \in [0,\infty)}|r_{kl}(t)|<1
\end{eqnarray}
\iE{will be further assumed}. For simplicity we consider only two uniform grids $\mathfrak{R}(\delta_1)$ and $\mathfrak{R}(\delta_2)$. Recall that $\delta_i,i=1,2$
depend on $T>0$; in the case of Pickands grid we set
$$\RDI= \mathfrak{R}(D_i a_T^{-2/\alpha})$$
 for some \iE{constant} $D_i>0, i=1,2$.  The vector of maxima on continuous time will be denoted by $\vk{M}(T)$ and that with respect to the discrete \AH{uniform} grid
$\mathfrak{R}(\delta_i),i=1,2$ \AH{by} $\vk{M}(\iH{\delta_i}, T)$. This means that the $k$th components of these two random vectors
are $M_k(T)$ and $M_k(\iH{\delta_i}, T)$, respectively which are  defined by
$$M_{k}(T)=\max_{t\in [0,T]} X_{k}(t), \quad M_{k}(\iH{\delta_i}, T)=\max_{t\in \RDI \cap[0,T]} X_{k}(t), \quad {k\le p}.$$
For notational simplicity we shall set below
$$ \widetilde{\vk{M}}(T)= \Bigl(a_T(M_1(T)- b_T) \ldot  a_T(M_p(T)- b_T)\Bigr)$$
and
$$\widetilde{\vk{M}}(\iH{\delta_i},T)=\Bigl(a_T(M_1(\iH{\delta_i}, T)- b_{\delta_i, T})\ldot
a_T(M_p(\iH{\delta_i}, T)- b_{\delta_i, T})\Bigr),$$
where $b_{\delta_i,T}$ is defined in \eqref{deltaT} if the grid $\mathfrak{R}(\delta_i)$ is sparse,
$b_{ \delta_i,T}= b_{T}(D_i)$ is given by \eqref{bDT}
if we consider a Pickands grid $\mathfrak{R}(\delta_i)=\mathfrak{R}(D_i a_T^{-2/\alpha})$ and
for a dense grid \bE{we set} $b_{\delta_i,T}=b_T$ with $b_T$ defined in \eqref{bT}. \\
In the following $\vk{x}, \vk{y}_1, \vk{y}_2\inr^p$ are fixed vectors and
$\vk{Z}$ is a $p$-dimensional centered Gaussian random vector with covariances
\BQN \label{covZ}
Cov(Z_{k},Z_{l})&=&\frac{r_{kl}}{\sqrt{r_{kk}r_{ll}}} , \quad 1 \le l \le k \le p.
\EQN
\bE{When $r_{kk}r_{ll}=0$ we assume that $Z_k$ and $Z_l$ are independent, i.e., we shall set
$$Cov(Z_{k},Z_{l})=0.$$}
The operations with vectors are meant componentwise, for instance $\x\le \vk{y}$ means $x_k \le y_k$ for any index $k\le p$, with $x_k$ and $y_k$
the $k$th component of $\x$ and $\y$, respectively. Hereafter we define
\def\pXY{p_{T,\x,\y, \delta}}
$$ \pXY:=\pk{ \widetilde{\M}(T)\le \x, \widetilde{\M}(\delta_i,T)\le \y_i, i=1,2}.$$

In the first \aE{theorem} below we \bE{discuss} the case when one of the grids is sparse. Our results shall establish that
\begin{eqnarray}
\label{eqAAA}
\lim_{T\to \infty}\sup_{ \xyrp } \ABs{
\pXY - \EE{\exp\Bigl(- \sum_{k=1}^{p}
f(x_k,y_{k1},y_{k2}) e^{-r_{kk}+\sqrt{2r_{kk}}Z_{k}}\Bigr)} }=0,
\end{eqnarray}
where the function $f$ is given below explicitly for each particular case.

\BT \label{Th1}
Let $\{\X(t), t\ge 0\}$ be \aE{a} centered stationary Gaussian vector process as defined above and let $\mathfrak{R(\delta_{1})}$ be a sparse grid.  Assume \iH{that}
(\ref{eq1.3}), (\ref{eq1.4}) and (\ref{eq1.5}) hold and the Gaussian random vector $\vk{Z}$ has a positive-definite covariance matrix \aE{with elements defined in \eqref{covZ}.} \\
i) If $\mathfrak{R(\delta_{2})}$  is another sparse grid such that $\mathfrak{R(\delta_{1})}\cap\mathfrak{R(\delta_{2})}=\emptyset$ or $\limit{T} \delta_1(T)/ \delta_2(T)=\infty$,
 then  \eqref{eqAAA} holds with
\BQNY
f(x_k,y_{k1}, y_{k2}) &=&
e^{-x_{k}}+e^{-y_{k1}} +e^{-y_{k2}}.
\EQNY
ii) 
\bE{Let} $\mathfrak{R(\delta_{2})}$  be a sparse grid such that $\mathfrak{R(\delta_{1})}\cap\mathfrak{R(\delta_{2})}
=\mathfrak{R(\delta_{3})}$. \bE{\iE{If}  $\mathfrak{R(\delta_{3})}$ is a \iE{non-empty} grid such that}
$$\lim_{T\to \IF} \ln(\frac{\delta_{3}(T)}{\delta_{1}(T)})= \theta_{1} \in [0,\IF), \quad
\lim_{T\to \IF} \ln(\frac{\delta_{3}(T)}{\delta_{2}(T)})= \theta_{2} \in [0,\IF),$$
then \eqref{eqAAA} holds with (write $\theta=\theta_{2}-\theta_{1})$
\BQNY
f(x_k,y_{k1}, y_{k2}) &=&e^{-x_k} + e^{-y_{k1}}+e^{-y_{k2}}-
e^{-y_{k1}- \aE{\theta_{1}}} I(y_{k1}>y_{k2}+\aE{\theta})-e^{-y_{k2}-\theta_{2}}I(y_{k1}\leq y_{k2}+\aE{\theta}),
\EQNY
where $I(\cdot)$ is the indicator function.

\TT{iii) If $\mathfrak{R(\delta_{2})}= \mathfrak{R}(D_2 a_T^{-2/\alpha})$ is a Pickands grid, then  \eqref{eqAAA} holds with
\BQNY
f(x_k,y_{k1},y_{k2})&=&
 e^{-x_k} + e^{-y_{k1}}+ e^{-y_{k2}}-H_{D_{2},\alpha}^{\ln H_{\alpha}+x_{k},\ln H_{D_{2},\alpha}+y_{k2}}.
\EQNY}

iv) If $\mathfrak{R(\delta_{2})}$ is a dense grid, then  again \eqref{eqAAA} holds with
\BQNY
\TT{f(x_k,y_{k1},y_{k2})=e^{-y_{k1}}+e^{-\min(x_k,y_{k2})}.}
\EQNY
\ET

We \aE{consider} next the \aE{cases} that one grid is \aE{a} Pickands grid \aE{and the \iE{second one} is either a Pickands or a dense grid.}
For positive constants $D_1,D_2,\lambda$  and $x, z_1,z_2  \aE{
\inr}$
define (recall $B^*_{\alpha/2}(t):= \sqrt{2}B_{\alpha/2}(t)-\iH{\abs{t}}^{\alpha}$)
\Ta{$$H_{D_1,D_2, \alpha}^{z_1,z_2}(\lambda)=
\int_{s\inr}e^{s} \pk{
\max_{k\inn: kD_i\in[0,\lambda]}B^*_{\alpha/2}(kD_i)>s+z_i,i=1,2
}\, ds$$}
and
$$H_{D_1,D_2, \alpha}^{x,z_1,z_2}(\lambda)=
\int_{s\inr}e^{s} \pk{ \max_{t\in [0,\lambda]} B^*_{\alpha/2}(t)> s+x,
\max_{k\inn: kD_i\in[0,\lambda]}B^*_{\alpha/2}(kD_i)>s+z_i,i=1,2
}\, ds.$$

\BT\label{Th2} Under the assumptions of Theorem \ref{Th1} \bE{suppose further}  that $\mathfrak{R(\delta_{1})}=\mathfrak{R}(D_1 a_T^{-2/\alpha}),D_1>0$ is a Pickands grid.

i) If $\mathfrak{R(\delta_{2})}=\mathfrak{R}(D_2 a_T ^{-2/\alpha}),\aE{D_2 \in (0,\IF) \setminus \{D_1\}}$ is also \iE{a} Pickands grid,
then for any $\aE{x,z_1,z_2\inr}$ 
$$\Ta{H_{D_1,D_2,\alpha}^{z_1,z_2}= \lim_{\lambda\to \IF}\frac{
H_{D_1,D_2,\alpha}^{z_1,z_2}(\lambda)}{\lambda}\in (0,\IF)\ \ \mbox{and}\ \ }
H_{D_1,D_2,\alpha}^{x,z_1,z_2}= \lim_{\lambda\to \IF}\frac{
H_{D_1,D_2,\alpha}^{x,z_1,z_2}(\lambda)}{\lambda}\in (0,\IF)$$
and further \eqref{eqAAA} holds with $f$ given by
\BQNY
f(x_k,y_{k1},y_{k2})&=&
e^{-x_k} + e^{-y_{k1}}+e^{-y_{k2}} - H_{D_1,\alpha}^{\ln H_\alpha+x_k,\ln H_{D_1,\alpha}+ y_{k1}} -
H_{D_2,\alpha}^{\ln H_\alpha+x_k,\ln H_{D_2,\alpha}+y_{k2} } \\
&&- H_{D_1,D_2,\alpha}^{
\ln H_{D_1,\alpha}+y_{k1},\ln H_{D_2,\alpha}+y_{k2} } +  H_{D_1,D_2,\alpha}^{\ln H_{\alpha}+x_k, \ln H_{D_1,\alpha}+y_{k1},\ln H_{D_2,\alpha}+ y_{k2}}.
\EQNY
\COM{
\BQN
\label{eq2.2}
\lim_{T\to \IF}\sup_{ \xyrp } \ABs{
P\left\{a_{T}(M_{k}^{\delta_{1}}(T)-b_{c,T})\leq x_{k}, a_{T}(M_{k}^{\delta_{2}}(T)- b_{d,T})\leq y_{k}, k=1,\cdots,p\right\}
-\EE{\exp\left(-g(\mathbf{x}, \mathbf{y}, \mathbf{Z})\right)} }=0,
\EQN
where
$$g(\mathbf{x}, \mathbf{y}, \mathbf{Z})=\sum_{k=1}^{p}\left(e^{-x_{k}-r_{kk}+\sqrt{2r_{kk}}Z_{k}}+e^{-y_{k}-r_{kk}+\sqrt{2r_{kk}}Z_{k}}-H_{c,d,\alpha}
^{ \ln   H_{c,\alpha}+x_{k}, \ln   H_{d,\alpha}+y_{k}}e^{-r_{kk}+\sqrt{2r_{kk}}Z_{k}}\right),$$
with
\BQN\label{btd2}
b_{c,T}=a_T+ a_T  \ln  \Bigl((2\pi)^{-1/2}C^{1/\alpha}H_{c,\alpha}a_T^{-1+2/\alpha}\Bigr).
\EQN
}
ii) If $\mathfrak{R(\delta_{2})}$ is a dense grid, then \eqref{eqAAA} holds with
\BQNY
f(x_k,y_{k1},y_{k2})&=&
\TT{e^{-\min(x_k,y_{k2})}+ e^{-y_{k1}} - H_{D_1,\alpha}^{\ln H_\alpha+\min(x_k,y_{k2}), \ln H_{D_1,\alpha}+ y_{k1}}.}
\EQNY
iii) If both $\mathfrak{R(\delta_{1})}$ and $\mathfrak{R(\delta_{2})}$  are dense \Quan{grids}, then again
\eqref{eqAAA} holds with
$$f(x_k,y_{k1},y_{k2})= e^{-\min(x_k,y_{k1},y_{k2})}.$$
\ET

\aE{{\bf Remarks}:  a) From the above results it follows that the joint convergence stated therein is determined by the choice of the grids.
The dependence parameters $r_{lk}, l,k\le p$ determine the covariance of the Gaussian random vector $\vk{Z}$ \bE{and appears explicitly in the
definition of the limiting distribution.} \\
Clearly, if each $r_{kk}$ equals 0, i.e., the Berman condition holds for each component of the vector process,
 then $\vk{Z}$ does not appear in any of the \jE{limiting results} above. \bE{For such cases} the maxima over a sparse grid is independent of
 that taken over a Pickands or a dense grid. \\
b)  Condition (9) can be stated in a slightly more general form putting \bE{therein} $C_k$ instead of $C$. Our results can be restated then
with some obvious \jE{modifications} on the constants involved. \\
c) In \cite{TanTang} a particular case of Piterbarg's max-discretisation theorem \bE{was \jE{investigated}}, which in our \AH{notation} corresponds to $r_{kk}=\IF$.
Considering for simplicity $p=1$, so we assume that $r_{11}=\IF$, then if \eqref{eq1.1} holds with $\alpha \in (0,1]$ \jE{and} $r(t)=o(1), t\to \IF$ a convex function, and $(r(t) \ln t)^{-1}$ is monotone for large $t$ and $o(1)$, then for any two \bE{different} sparse, Pickands or dense grids $\mathfrak{R(\delta_{1})}$ and $\mathfrak{R(\delta_{2})}$  we have
\BQN
\quad \iH{\limit{T}}\pk{a_T^*(M(T) - b_{T}^*) \le x, a_T^*(M(\iH{\delta_1},T) - b_{\delta_1, T}^*) \le y, a_T^*(M(\iH{\delta_2},T) - b_{\delta_2,T}^*) \le z}
= \Phi( \min(x,y,z))
\EQN
for any $x,y,z\inr$ as $T\to \IF$, where
$$a_T^*=1/\sqrt{r(T)}, \quad b_{\delta_i,T}^*=\sqrt{ (1 - r(T))/r(T)} b_{ \delta_i, T}$$
 and $\Phi$ denotes the distribution function of an
 $N(0,1)$ random variable. The proof of the above claim follows by Theorem 2.1 in \cite{TanTang} and \nelem{lem:A4}. \\
\bE{Consequently, for this case different grids do not play a role in the limiting distribution. Note however that the noramlisation constant $b_{\delta_i,T}^*$ depends
on the \jE{type of the} grid.}}\\
\aE{iv) Set for $\x,\y_1, \y_2 \inr^p$
$$ G(\x,\y_1,\y_2)= \EE{\exp\Bigl(- \sum_{k=1}^{p}
f(x_k,y_{k1},y_{k2}) e^{-r_{kk}+\sqrt{2r_{kk}}Z_{k}}\Bigr)} ,$$
 where $f$ and $\vk{Z}$ are \iE{as} in Theorem 2.1 and Theorem 2.2. It follows that
\BQNY
\lim_{x\to - \IF} H_{D_1,D_2,\alpha}^{x,y_1,y_2}= H_{D_1,D_2,\alpha}^{y_1,y_2}, \quad
\lim_{y_1 \to - \IF, y_2 \to - \IF} H_{D_1,D_2,\alpha}^{x,y,z }= e^{-x} H_\alpha,\\
\lim_{x \to - \IF, y_1 \to - \IF} H_{D_1,D_2,\alpha}^{x,y_1,y_2}=
\lim_{y_1 \to - \IF} H_{D_1,D_2,\alpha}^{y_1, y_2}= H_{D_2,\alpha}^{y_2}= e^{-y_2} H_{D_2,\alpha}.
\EQNY
}
Hence, using further \eqref{eqzIF} we conclude that $G$ is a non-degenerate multivariate distribution in $\R^{3p}$, which we refer to as the Piterbarg distribution. One important property of \bE{$G$ is} that when $r_{kk}=0$ for all indices $k\le p$, then \bE{it has} unit Gumbel marginals $\Lambda(x)= e^{-e^{-x}},x\inr$. Moreover, \bE{$G$ is} a max-stable \jE{distribution} since
$$ (G(\x+ \ln n, \y_1 + \ln n, \y_2 + \ln n))^n= G(\x,\y_1,\y_2), \quad \x_1,\y_1,\y_2 \inr^p, n\inn.$$
In Extreme Value Theory max-stable \jE{distributions} are important for modelling of extremes and rare events, see e.g., \cite{Res1987,Faletal2010}
for details.

\section{Proofs}

In this section we present several lemmas needed for the proof
of the main results. In order to establish Piterbarg's max-discretisation theorem for multivariate stationary Gaussian processes we need to closely follow
\cite{Pit2004}, and of course to
strongly rely on the deep ideas and the techniques presented in \cite{Pit96}. First, for $1\leq k,l\leq p$ define
$$\rho_{kl}(T)=r_{kl}/\ln T.$$
Following the former reference, we divide the interval $[0,T]$ onto
intervals of length $S$ alternating with shorter intervals
of length $R$. Let \AH{$ b< a < 1 $} be \aE{two positive} constants, \aE{where $b$ will be \bE{chosen} below (\iE{see}  \eqref{bkk})}.
We shall denote throughout in the sequel
$$S=T^{a}, \quad  R=T^{b}, \quad T>0.$$
 Denote the long intervals by $ \SL $, $l=1,\cdots,n_{T}$, and the short intervals by $ \RL $,
$l=1,\cdots,n_{T}$ where
\BQN \label{nT}
n_T:=[T/(S+R)].
\EQN
  It will be seen from the proofs, that a
possible remaining interval with length different than $S$ or $R$
plays no role in our asymptotic considerations; we call also
this interval a short interval. Define further $\mathbf{S}=\cup_{l=1}^{n_T}
\SL , \mathbf{R}=\cup_{l=1}^{n_T}  \RL $ and thus $[0,T]=\mathbf{S}\cup
\mathbf{R}$.

Our proofs also \bE{rely}  on the ideas of \cite{mittal1975limit};  \aE{we shall  construct}
  new Gaussian processes to approximate the
original ones. For each index $k \le p$ we define a Gaussian process $\eta_k$ as
\BQN \label{etak}
\eta_{k}(t)=Y_{k}^{(j)}(t), \quad t\in \RJ\cup\SJ=[(j-1)(S+R),j(S+R)),
\EQN
where  $\{Y_{k}^{(j)}(t), t\geq
0\}$, $j=1,\cdots, n_{T}$ are independent copies of $\{X_{k}(t),
t\geq0\}$. We construct the processes so that
$\eta_k,k=1, \cdots, p$ are independent by taking $Y_k^{(j)}$ to
be independent for any $j$ and $k$ two possible indices. The
independence of  $\eta_k$ and $\eta_l$
implies
$$\gamma_{kl}(s,t):=\EE{\eta_{k}(s)\eta_{l}(t)}=0,\ \ k\neq l,$$
whereas for any fixed $k$
\BQNY
  \gamma_{kk}(s,t)&:=&\EE{\eta_{k}(s)\eta_{k}(t)}\\
  &=&\left\{
 \begin{array}{cc}
  \EE{ Y_k^{(i)}(t),Y_k^{(i)}(s)}= { r_{kk}(s,t)},   & \text{ if } t,s\in \mathcal{R}_{i}\cup \mathcal{S}_{i}, \text{ for some } i \le n_{T};\\
  \EE{Y_k^{(i)}(t),Y_k^{(j)}(s)}={0},    & \text{ if } t\in \mathcal{R}_{i}\cup \mathcal{S}_{i}, s\in \mathcal{R}_{j}\cup \mathcal{S}_{j},
  \text{ for some } i\neq j \le n_{T}.
 \end{array}
  \right.
\EQNY
For $k=1,\cdots,p$ define
$$\xi_{k}^{T}(t)=\big(1-\rho_{kk}(T)\big)^{1/2}\eta_{k}(t)+\rho^{1/2}_{kk}(T)Z_{k}, \ \ 0\leq t\leq T,$$
where $\mathbf{Z}=(Z_1, \ldots, Z_p)$ is a $p$-dimensional
centered Gaussian random vector \jE{introduced} in Section
2, which is independent of $\{\eta_{k}(t), t\geq 0\}$,
$k=1,\cdots,p$. \COM{It can be checked that  $\{\xi_{k}^{T}(t),
0\leq t\leq T, k=1,\cdots,p\}$ are stationary Gaussian vector
processes.} Denote by $\{\varrho_{kl}(s,t),1\leq k,l\leq p\}$ the
covariance functions of $\{\xi_{k}^{T}(t), 0\leq t\leq
T,k=1,\cdots,p\}$. We have
$$\varrho_{kl}(s,t)=\EE{\xi_{k}^{T}(s)\xi_{l}^{T}(t)}=\rho_{kl}(T),\ \ k\neq l$$
and
\[
  \varrho_{kk}(s,t)=\left\{
 \begin{array}{cc}
  {r_{kk}(s,t)+(1-r_{kk}(s,t))\rho_{kk}(T)} ,    &t\in \mathcal{R}_{i}\cup \mathcal{S}_{i}, s\in \mathcal{R}_{j}\cup \mathcal{S}_{j}, i= j;\\
  {\rho_{kk}(T)},    & t\in \mathcal{R}_{i}\cup \mathcal{S}_{i}, s\in \mathcal{R}_{j}\cup \mathcal{S}_{j}, i\neq j.
 \end{array}
  \right.
\]
\COM{Since both $\gamma_{kk}(s,t)$ and $\varrho_{kk}(s,t)$ are
functions of $|s-t|$  we shall write below simply
$\gamma_{kk}(|s-t|)$ and $\varrho_{kk}(|s-t|)$, respectively.}

For any $\varepsilon>0$ set
\BQN \label{qE}
\LE{q_\ve}= \frac{\varepsilon}{ ( \ln   T)^{1/\alpha}}.
\EQN
For notational simplicity we write
$$\widetilde{\vk{M}}_{\xi}(\LE{q_\ve},\mathbf{S})=\Bigl(a_T(M_{\xi 1}(\LE{q_\ve},\mathbf{S} )- b_{T})\ldot
a_T(M_{\xi p}(\LE{q_\ve},\mathbf{S})- b_{T})\Bigr)$$
and
$$\widetilde{\vk{M}}_{\xi}(\delta_i,\mathbf{S})=\Bigl(a_T(M_{\xi 1}(\delta_i,\mathbf{S} )- b_{\delta_i,T})\ldot
a_T(M_{\xi p}(\delta_i,\mathbf{S})- b_{\delta_i,T})\Bigr),$$
where $$M_{\xi k}(\LE{q_\ve}, \mathbf{S})=\max_{t\in \mathfrak{R}(\iE{q_\ve})\cap \mathbf{S}}\xi_{k}^{T}(t)$$
and $b_{\delta_i,T}$ is defined in \eqref{deltaT} if the grid $\mathfrak{R}(\delta_i)$ is sparse,
$b_{ \delta_i,T}= b_{T}(D_i)$ is given by \eqref{bDT}
if we consider a Pickands grid $\mathfrak{R}(\delta_i)=\mathfrak{R}(D_i a_T^{-2/\alpha})$ and
for a dense grid $b_{\delta_i,T}=b_T$ with $b_T$ defined in \eqref{bT}.

We \aE{present} first four lemmas. \aE{Since their proofs are similar to those of}  Lemmas 3.1-3.4 in \cite{TanH2014} we shall
\aE{not give them here.}

\BL\label{lem:3.1}  \jE{If} $\mathfrak{R}(\delta_{1})$ and $\mathfrak{R}(\delta_{2})$ are sparse   or Pickands \aE{grids}, then for any
$B>0$ there exits some $K>0$ such that for
all $x_{k},y_{ki}\in[-B,B], i=1,\jE{2},  k\le p$
\begin{eqnarray*}
&&\bigg|\pk{ \widetilde{\M}(T)\le \x, \widetilde{\M}(\delta_i,T)\le \y_i, i=1,2}-\pk{ \widetilde{\M}(\mathbf{S})\le \x, \widetilde{\M}(\delta_i,\mathbf{S})\le \y_i, i=1,2}\bigg| \leq K( \ln   T)^{1/\alpha-1/2}T^{b-a}
\end{eqnarray*}
holds for some $0< b< a < 1$ and all $T$ large.
\EL

\LE{In the following $\LE{\mathfrak{R}(q_\ve)}=\mathfrak{R}(\varepsilon/(\ln T)^{1/\alpha})$ denotes a Pickands grid where $\ve>0$ \iE{and $q_\ve$ is defined in \eqref{qE}}.}

\BL\label{lem:3.2} \jE{If}
$\mathfrak{R}(\delta_{1})$ and $\mathfrak{R}(\delta_{2})$ are sparse or Pickands \aE{grids}, then for any
$B>0$ \Quan{and} for all $x_{k},y_{ki}\in[-B,B], i=1,2,  k\le p$
\begin{eqnarray*}
&&\bigg|\pk{ \widetilde{\M}(\mathbf{S})\le \x, \widetilde{\M}(\delta_i,\mathbf{S})\le \y_i, i=1,2}-\pk{ \widetilde{\M}(q_\ve,\mathbf{S})\le \x, \widetilde{\M}(\delta_i,\mathbf{S})\le \y_i, i=1,2}\bigg|\rightarrow 0
\end{eqnarray*}
as $\varepsilon\downarrow 0$.
\EL

\BL\label{lem:3.3} \jE{If}
$\mathfrak{R}(\delta_{1})$ and $\mathfrak{R}(\delta_{2})$ are sparse or \aE{Pickands grids}, then  for any
$B>0$ \Quan{and} for all $x_{k},y_{ki}\in[-B,B], i=1,2,  k\le p$
\begin{eqnarray*}
&& \lim_{T \to \IF} \bigg|\pk{ \widetilde{\M}(\LE{q_\ve},\mathbf{S})\le \x, \widetilde{\M}(\delta_i,\mathbf{S})\le \y_i, i=1,2}
-\pk{ \widetilde{\M}_{\xi}(\LE{q_\ve},\mathbf{S})\le \x, \widetilde{\M}_{\xi}(\delta_i,\mathbf{S})\le \y_i, i=1,2}\bigg|
=0
\end{eqnarray*}
\Tan{uniformly for $\varepsilon>0$.}
\EL

\jE{Let in the following $\Phi_p$ denote the distribution \Quan{function
of}  the $p$-dimensional Gaussian random vector $\mathbf{Z}$ and set for} $\eta_k$ defined in \eqref{etak}
$$\widehat{\M}_{\eta}(\delta_i,\mathcal{S}_{j})=\left(\max_{t\in \RDI\cap \mathcal{S}_{j}}\eta_{1}(t),\cdots, \max_{t\in \RDI\cap \mathcal{S}_{j}}\eta_{p}(t)\right), \quad \widehat{\M}_{\eta}(\mathcal{S}_{j})=\left(\max_{t\in \mathcal{S}_{j}}\eta_{1}(t),\cdots, \max_{t\in \mathcal{S}_{j}}\eta_{p}(t)\right).$$

\BL\label{lem:3.4} \jE{If}
$\mathfrak{R}(\delta_{1})$ and $\mathfrak{R}(\delta_{2})$ are sparse or Pickands \aE{grids}, then for
any $B>0$ for all $x_{k},y_{ki}\in[-B,B], i=1,2,  k\le p$
\begin{eqnarray*}
&&\bigg|
\pk{ \widetilde{\M}_{\xi}(\LE{q_\ve},\mathbf{S})\le \x, \widetilde{\M}_{\xi}(\delta_i,\mathbf{S})\le \y_i, i=1,2}
\\
&&\ \ \ \ \ \ - \int_{\mathbf{z} \in \mathbb{R}^p}\prod_{j=1}^{n_{T}}\pk{ \widehat{\M}_{\eta}(\mathcal{S}_{j})\le \mathbf{u}(\x,\vk{z}), \widehat{\M}_{\eta}(\delta_i,\mathcal{S}_{j})\le \mathbf{u}(\y_i,\vk{z}), i=1,2}
d\Phi_{p}(\vk{z})
\bigg|\rightarrow 0
\end{eqnarray*}
as $\varepsilon\downarrow 0$, where  $\vk{u}(\x,\vk{z}), \vk{u}(\y_i,\vk{z}),i=1,2$ have components
\begin{eqnarray}
\label{eq3.4.2}
u(x_{k},z_k)=\frac{b_{T}+x_{k}/a_{T}-\rho^{1/2}_{kk}(T)z_{k}}{(1-\rho_{kk}(T))^{1/2}}=\frac{x_{k}+r_{kk}-\sqrt{2r_{kk}}z_{k}}{a_{T}}+b_{T}+o(a_{T}^{-1}),
\end{eqnarray}
\begin{eqnarray}
\label{eq3.4.3}
u(y_{ki}, \LE{z_k})=\frac{b_{\delta_i,T}+y_{ki}/a_{T}-\rho^{1/2}_{kk}(T)z_{k}}{(1-\rho_{kk}(T))^{1/2}}=\frac{y_{ki}+r_{kk}-\sqrt{2r_{kk}}z_{k}}{a_{T}}+b_{\delta_i,T}+o(a_{T}^{-1}),
\end{eqnarray}
\Tan{for all $x_{k},y_{ki}\in[-B,B], i=1,2,  k\le p$.}
\EL

\def\AK{\mathcal{A}_{k}}
\def\AL{\mathcal{A}_{l}}
\def\AJ{\mathcal{A}_{j}}
\def\APK{\mathcal{A}_{p+k}}
\def\APL{\mathcal{A}_{p+l}}
\def\APJ{\mathcal{A}_{p+j}}
\def\APPK{\mathcal{A}_{2p+k}}
\def\APPL{\mathcal{A}_{2p+l}}
\def\APPJ{\mathcal{A}_{2p+j}}

\prooftheo{Th1} \Tan{Since all the limits of the probabilities in Lemmas 3.1-3.4 \jE{are positive} for all $x_{k},y_{ki}\in[-B,B], i=1,2,  k\le p$, by letting $\varepsilon\downarrow 0$, we have
\begin{eqnarray*}
&&\pk{ \widetilde{\M}(T)\le \x, \widetilde{\M}(\delta_i,T)\le \y_i, i=1,2}\\
&&\ \ \sim\int_{\mathbf{z} \in \mathbb{R}^p}\prod_{j=1}^{n_{T}}\pk{ \widehat{\M}_{\eta}(\mathcal{S}_{j})\le \mathbf{u}(\x,\vk{z}), \widehat{\M}_{\eta}(\delta_i,\mathcal{S}_{j})\le \mathbf{u}(\y_i,\vk{z}), i=1,2}
d\Phi_{p}(\vk{z})
\end{eqnarray*}
as $T\rightarrow\infty$. Thus, if we can prove
\begin{eqnarray}
\label{eq2.1.1}
&&\lim_{T\rightarrow\infty}\bigg|\prod_{j=1}^{n_{T}}\pk{ \widehat{\M}_{\eta}(\mathcal{S}_{j})\le \mathbf{u}(\x,\vk{z}), \widehat{\M}_{\eta}(\delta_i,\mathcal{S}_{j})\le \mathbf{u}(\y_i,\vk{z}), i=1,2} \nonumber\\
&&\ \ \ \ \ \ \ \ \ \ \ \ \ \ \ \ \ \ \ \ \ -\exp\Bigl(- \sum_{k=1}^{p}
f(x_k,y_{k1},y_{k2}) e^{-r_{kk}+\sqrt{2r_{kk}}z_{k}}\Bigr)\bigg|=0,
\end{eqnarray}
where $f(x_k,y_{k1},y_{k2})$ is defined in Theorem 2.1, then applying the dominated
convergence theorem we complete the proof of Theorem 2.1 for the case $i)-iii)$.}
Define next the events
$${\mathcal{A}}_k=\Bigl\{ \max_{t\in
[0,S]}\eta_{k}(t)> u(x_{k},z_{k}) \Bigr\}, \quad  {\mathcal{A}}_{p+k}=
\Bigl\{\max_{t\in\mathfrak{R}(\delta_{1})\cap [0,S]}\eta_{k}(t)>
u(y_{k1},z_{k})\Bigr\}$$ and $$ {\mathcal{A}}_{2p+k}=
\Bigl\{\max_{t\in\mathfrak{R}(\delta_{2})\cap [0,S]}\eta_{k}(t)>
u(y_{k2},z_{k})\Bigr\},\quad  k=1,\cdots,p.$$

$i)$ Using the stationarity of
$\{\eta_{k}(t), k=1,\cdots,p\}$ (we write $\AK^c$ for the
complimentary event of $\AK$)
\BQNY
\prod_{j=1}^{n_{T}}\pk{\widehat{\M}_{\eta}(\mathcal{S}_{j})\le \mathbf{u}(\x,\vk{z}), \widehat{\M}_{\eta}(\delta_i,\mathcal{S}_{j})\le
 \mathbf{u}(\y_i,\vk{z}), i=1,2}
 &=& (\mathbb{P}\{ \cap_{k=1}^{3p} \AK^c \})^{n_{T}}\\
&=&\exp\big(n_{T} \ln  (  \mathbb{P}\{ \cap_{k=1}^{3p} \AK^c \})\big)\\
&=&\exp\big(-n_{T} \mathbb{P}\{ \cup_{k=1}^{3p} \AK \}+W_{n_{T}}\big),
 \EQNY
 where $n_T$ is defined in \eqref{nT}.
Since $\lim_{T \to \IF} \mathbb{P}\{ \cap_{k=1}^{3p} \AK \} =1$ we get that
the remainder $W_{n_{T}}$ satisfies
$$W_{n_{T}}=o(n_{T} \mathbb{P}\{ \cup_{k=1}^{3p} \AK \}), \quad T \to \IF.$$
Next, by  Bonferroni inequality
\BQN\label{bonfi}
\sum_{k=1}^{3p} \mathbb{P}\{ \AK \}  &\ge & \mathbb{P}\{ \cup_{k=1}^{3p} \AK \} 
 \ge  \sum_{k=1}^{3p} \mathbb{P}\{ \AK \} - \sum_{1 \le k < l \le 3p} \mathbb{P}\{ \AK, \AL \} \notag\\
& = & \sum_{k=1}^{3p} \mathbb{P}\{ \AK \} - \sum_{1 \le k < l \le p} \mathbb{P}\{ \AK, \AL \} - \sum_{1 \le k < l \le p} \mathbb{P}\{ \APK, \APL \} -\sum_{1 \le k < l \le p} \mathbb{P}\{ \APPK, \APPL \}\notag\\
 &&- 2 \sum_{1 \le k < l \le p} \mathbb{P}\{ \AK, \APL \}- 2 \sum_{1 \le k < l \le p} \mathbb{P}\{ \AK, \APPL \} -2 \sum_{1 \le k < l \le p} \mathbb{P}\{ \APK, \APPL \}\notag\\
&&- \sum_{k=1}^p  \mathbb{P}\{ \AK,\APK\}-\sum_{k=1}^p  \mathbb{P}\{ \AK,\APPK\}-\sum_{k=1}^p  \mathbb{P}\{ \APK,\APPK\}\notag\\
&=:&A_{1}-A_{2}-A_{3}-A_{4}-2A_{5}-2A_{6}-2A_{7}-A_{8}-A_{9}-A_{10}. \EQN

Further, Lemma 2 in  \cite{Pit2004} and (\ref{eq3.4.2}),
(\ref{eq3.4.3}) imply \AH{(recall $S=T^a$)}
\begin{eqnarray*}
\cE{A_{1}}&\sim & \sum_{k=1}^{p}ST^{-1}(e^{-x_{k}}+e^{-y_{k1}}+e^{-y_{k2}})e^{-r_{kk}+\sqrt{2r_{kk}}z_{k}}
, \quad T\rightarrow\infty.
\end{eqnarray*}
For $A_{2}$, by the independence of $\eta_{k}(t)$ and $\eta_{l}(t)$,
$k\neq l$, Lemma 2 of  \cite{Pit2004} and (\ref{eq3.4.2}),
(\ref{eq3.4.3}),
 we have
\begin{eqnarray*}
A_{2}&=&\sum_{1\leq k<l\leq p}\mathbb{P}\left\{\max_{t\in [0,S]}\eta_{k}(t)> u(x_{k},z_{k}), \max_{t\in [0,S]}\eta_{l}(t)> u(x_{l},z_{l})\right\}\nonumber\\
&=&\sum_{1\leq k<l\leq p}\mathbb{P}\left\{\max_{t\in [0,S]}\eta_{k}(t)> u(x_{k},z_{k})\right\}
\mathbb{P}\left\{ \max_{t\in [0,S]}\eta_{l}(t)> u(x_{l},z_{l})\right\}\nonumber\\
&\sim&\sum_{1\leq k<l\leq p}ST^{-1}e^{-x_{k}-r_{kk}+\sqrt{2r_{kk}}z_{k}}ST^{-1}e^{-y_{l}-r_{ll}+\sqrt{2r_{ll}}z_{l}}
=o(A_{1}).
\end{eqnarray*}
\LE{Since} $\mathfrak{R}(\delta_{i}), i=1,2$ is a sparse grid, similar arguments
as for $A_{2}$ lead to
$$A_{k}=o(A_{1}), \ \ k=3,4,5,6,7.$$
Further, Lemma 2 of  \cite{Pit2004} implies $A_{i}=o(A_{1}), i=8,9.$
By the first assertion of \nelem{lem:A1} we have
\begin{eqnarray*}
A_{10}=\Quan{o(T^{a-1})}=o(A_{1}).
\end{eqnarray*}
Consequently, as $T\to \IF$
$$n_{T} \mathbb{P}\{ \cup_{k=1}^{3p} \AK \} \sim \sum_{k=1}^{p}
(e^{-x_{k}}+e^{-y_{k1}}+e^{-y_{k2}})e^{-r_{kk}+\sqrt{2r_{kk}}z_{k}},$$
which completes the proof of (\ref{eq2.1.1}).

$ii)$ We proceed as for the proof of case $i)$ using the lower bound \eqref{bonfi}; we have thus
\BQN\label{case2}
\mathbb{P}\{ \cup_{k=1}^{3p} \AK \}
&=&  \sum_{k=1}^{3p} \mathbb{P}\{ \AK \} - \sum_{1 \le k , l \le 3p} \mathbb{P}\{ \AK, \AL \}+ \sum_{1 \le k,\cdots, l \le 3p}\Tan{\mathbb{P}\{\cdot\}} \notag\\
& = & \sum_{k=1}^{3p} \mathbb{P}\{ \AK \} - \sum_{1 \le k < l \le p} \mathbb{P}\{ \AK, \AL \} - \sum_{1 \le k < l \le p} \mathbb{P}\{ \APK, \APL \} -\sum_{1 \le k < l \le p} \mathbb{P}\{ \APPK, \APPL \}\notag\\
 &&- 2 \sum_{1 \le k < l \le p} \mathbb{P}\{ \AK, \APL \}- 2 \sum_{1 \le k < l \le p} \mathbb{P}\{ \AK, \APPL \} -2 \sum_{1 \le k < l \le p} \mathbb{P}\{ \APK, \APPL \}\notag\\
&&- \sum_{k=1}^p  \mathbb{P}\{ \AK,\APK\}-\sum_{k=1}^p  \mathbb{P}\{ \AK,\APPK\}-\sum_{k=1}^p
 \mathbb{P}\{ \APK,\APPK\}+\sum_{1\leq k<,\cdots,l\leq 3p}\Tan{\pk{\cdot}} \notag\\
&=:&A_{1}-A_{2}-A_{3}-A_{4}-2A_{5}-2A_{6}-2A_{7}-A_{8}-A_{9}-A_{10}+A_{11}. \EQN
The estimates for $A_{i}$, $i=1,\cdots,9$ are the same as for case $i)$, therefore we only need to deal with the terms $A_{10}$ and $A_{11}$.
\jE{It follows that each term of} $A_{11}$ can be bounded by $A_{5}$, $A_{6}$ or $A_{7}$ implying
$$A_{11}=o(A_{1}), \quad T\to \IF.$$
Next, the definition of $u(y_{ki},z_{k}), i=1,2$ implies
\begin{eqnarray}
\label{eq.pr1}
u(y_{ki},z_{k})=\sqrt{2\ln T}-\frac{1}{2}\frac{\ln\ln T}{\sqrt{2\ln T}}+\frac{\ln \delta_{i}^{-1}(T)}{\sqrt{2\ln T}}+\frac{\ln (2^{-1}\pi^{-1/2})}{\sqrt{2\ln T}}+\Tan{\frac{y_{ki}+r_{kk}-\sqrt{2r_{kk}}z_{k}}{\sqrt{2\ln T}}+\frac{o_{T}(1)}{\sqrt{2\ln T}}}
\end{eqnarray}
for sparse grids. From the assumptions we know that $\lim_{T\to \IF} \ln(\frac{\delta_{1}\iH{(T)}}{\delta_{2}\iH{(T)}})= \theta=\theta_{2}-\theta_{1}$. Consequently, we have
\begin{eqnarray*}
u(y_{k1},z_{k})-u(y_{k2},z_{k})&=&\left[\ln(\frac{\delta_{2}(T)}{\delta_{1}(T)})+y_{k1}-y_{k2}\right](2\ln T)^{-1/2}\Tan{+o_{T}(1)(2\ln T)^{-1/2}}\\
&\sim&\left[-\theta +y_{k1}-y_{k2}\right](2\ln T)^{-1/2}\Tan{+o_{T}(1)(2\ln T)^{-1/2}}
\end{eqnarray*}
as $T\rightarrow\infty$.
Letting first $y_{k1}> y_{k2}+\theta$, we thus have $u(y_{k1},z_{k})> u(y_{k2},z_{k})$ for \jE{sufficiently} large $T$. Further,
\begin{eqnarray*}
A_{10}&=&\sum_{k=1}^{p}\mathbb{P}\left\{\max_{t\in \mathcal{R}(\delta_{1})\cap[0,S]}\eta_{k}(t)> u(y_{k1},z_{k}), \max_{t\in \mathcal{R}(\delta_{2})\cap[0,S]}\eta_{k}(t)> u(y_{k2},z_{k})\right\}\nonumber\\
&=&\sum_{k=1}^{p}\bigg[\mathbb{P}\left\{\max_{t\in \mathcal{R}(\delta_{1})\cap[0,S]}\eta_{k}(t)> u(y_{k1},z_{k})\right\}
+\mathbb{P}\left\{\max_{t\in \mathcal{R}(\delta_{2})\cap[0,S]}\eta_{k}(t)> u(y_{k2},z_{k})\right\}\\
&&-\left(1-\mathbb{P}\left\{\max_{t\in \mathcal{R}(\delta_{1})\cap[0,S]}\eta_{k}(t)\leq u(y_{k1},z_{k}), \max_{t\in \mathcal{R}(\delta_{2})\cap[0,S]}\eta_{k}(t)\leq u(y_{k2},z_{k})\right\}\right)\bigg]\\
&=&\sum_{k=1}^{p}\bigg[\mathbb{P}\left\{\max_{t\in \mathcal{R}(\delta_{1})\cap[0,S]}\eta_{k}(t)> u(y_{k1},z_{k})\right\}
+\mathbb{P}\left\{\max_{t\in \mathcal{R}(\delta_{2})\cap[0,S]}\eta_{k}(t)> u(y_{k2},z_{k})\right\}\\
&&-\left(1-\mathbb{P}\left\{\max_{t\in \mathcal{R}(\delta_{1})\cap[0,S]\setminus \mathcal{R}(\delta_{2})\cap[0,S]}\eta_{k}(t)\leq u(y_{k1},z_{k}),\max_{t\in \mathcal{R}(\delta_{2})\cap[0,S]}\eta_{k}(t)\leq u(y_{k2},z_{k})\right\}\right)\bigg]\\
&=&\sum_{k=1}^{p}\bigg[\mathbb{P}\left\{\max_{t\in \mathcal{R}(\delta_{1})\cap[0,S]\setminus \mathcal{R}(\delta_{3})\cap[0,S]}\eta_{k}(t)> u(y_{k1},z_{k})\right\}
-\mathbb{P}\left\{\max_{t\in \mathcal{R}(\delta_{1})\cap[0,S]}\eta_{k}(t)> u(y_{k1},z_{k})\right\}\\
&&+ \mathbb{P}\left\{\max_{t\in \mathcal{R}(\delta_{1})\cap[0,S]\setminus \mathcal{R}(\delta_{3})\cap[0,S]}\eta_{k}(t)> u(y_{k1},z_{k}), \max_{t\in \mathcal{R}(\delta_{2})\cap[0,S]}\eta_{k}(t)> u(y_{k2},z_{k})\right\}\bigg].
\end{eqnarray*}
\Quan{By \nelem{lem:A5} and (\ref{eq.pr1})
 we have for $i=1,2$ as $T\rightarrow\infty$
\begin{eqnarray*}
\mathbb{P}\left\{\max_{t\in \mathcal{R}(\delta_{\AH{\delta_i}})\cap[0,S]\setminus \mathcal{R}(\delta_{3})\cap[0,S]}\eta_{k}(t)> u(y_{k1},z_{k})\right\}
&\sim & \AH{ \delta _i \frac{S}{T}(\frac{1}{\AH{\delta_i} } - \frac{1}{\delta_3 }) e^{-y_{k1}-r_{kk}+\sqrt{2r_{kk}}z_{k}}}\\
&\sim &  ST^{-1} (1-\aE{e^{-\theta_{i}}})e^{-y_{k1}-r_{kk}+\sqrt{2r_{kk}}z_{k}}, \quad T\to \IF.
\end{eqnarray*}}
Further, applying Lemma 2 in \cite{Pit2004} (recall (\ref{eq.pr1})) we obtain as $T\rightarrow\infty$
\begin{eqnarray*}
\mathbb{P}\left\{\max_{t\in \mathcal{R}(\delta_{i})\cap[0,S]}\eta_{k}(t)> u(y_{ki},z_{k})\right\}\sim ST^{-1}e^{-y_{ki}-r_{kk}+\sqrt{2r_{kk}}z_{k}}, \quad \AH{i=1,2.}
\end{eqnarray*}
By the second assertion of \nelem{lem:A1}, the third term is $o(T^{a-1})$.\\
Next, for $y_{k1}\leq y_{k2}+ \theta$, we have $u(y_{k1},z_{k})\leq u(y_{k2},z_{k})$ for sufficient large $T$. Similarly, we have
\begin{eqnarray*}
A_{10}&=&\sum_{k=1}^{p}\mathbb{P}\left\{\max_{t\in \mathcal{R}(\delta_{1})\cap[0,S]}\eta_{k}(t)> u(y_{k1},z_{k}), \max_{t\in \mathcal{R}(\delta_{2})\cap[0,S]}\eta_{k}(t)> u(y_{k2},z_{k})\right\}\nonumber\\
&=&\sum_{k=1}^{p}\bigg[\mathbb{P}\left\{\max_{t\in \mathcal{R}(\delta_{2})\cap[0,S]}\eta_{k}(t)> u(y_{k2},z_{k})\right\}
-\mathbb{P}\left\{\max_{t\in \mathcal{R}(\delta_{2})\cap[0,S]\setminus \mathcal{R}(\delta_{3})\cap[0,S]}\eta_{k}(t)> u(y_{k2},z_{k})\right\}\\
&&+ \mathbb{P}\left\{\max_{t\in \mathcal{R}(\delta_{1})\cap[0,S]}\eta_{k}(t)> u(y_{k1},z_{k}), \max_{t\in \mathcal{R}(\delta_{2})\cap[0,S]\setminus \mathcal{R}(\delta_{3})\cap[0,S]}\eta_{k}(t)> u(y_{k2},z_{k})\right\}\bigg].
\end{eqnarray*}
\COM{\Quan{By \nelem{lem:A5} and (\ref{eq.pr1}) again
 we have as $T\rightarrow\infty$
 \begin{eqnarray*}
\mathbb{P}\left\{\max_{t\in \mathcal{R}(\delta_{2})\cap[0,S]\setminus \mathcal{R}(\delta_{3})\cap[0,S]}\eta_{k}(t)> u(y_{k2},z_{k})\right\}
&\sim& ST^{-1}(1-\aE{e^{-\theta_{2}}})e^{-y_{k2}-r_{kk}+\sqrt{2r_{kk}}z_{k}}
\end{eqnarray*}}
By Lemma 2 in \cite{Pit2004} and (\ref{eq.pr1}) again  we have as $T\rightarrow\infty$
\begin{eqnarray*}
\mathbb{P}\left\{\max_{t\in \mathcal{R}(\delta_{2})\cap[0,S]}\eta_{k}(t)> u(y_{k2},z_{k})\right\}\sim ST^{-1}e^{-y_{k2}-r_{kk}+\sqrt{2r_{kk}}z_{k}}.
\end{eqnarray*}
}
Again, in view of the second assertion of \nelem{lem:A1}  the third term is also $o(T^{a-1})$.
Consequently,
\begin{eqnarray*}
A_{10}&=&\sum_{k=1}^{p}T^{a-1}[e^{-y_{k1}- \aE{\theta_{1}}}I(y_{k1}> y_{k2}+ \theta)+e^{-y_{k2}- \aE{\theta_{2}}}I(y_{k1}\leq y_{k2}+ \theta)]e^{-r_{kk}+\sqrt{2r_{kk}}z_{k}}+o(\LE{T^{a-1}}), \quad T\rightarrow\infty
\end{eqnarray*}
implying that as $T\to \IF$
$$n_{T} \mathbb{P}\{ \cup_{k=1}^{3p} \AK \} \sim \sum_{k=1}^{p}(e^{-x_k} + e^{-y_{k1}}+e^{-y_{k2}}-
e^{-y_{k1}- \aE{\theta_{1}}} I(y_{k1}>y_{k2}+\aE{\theta})-e^{-y_{k2}-\theta_{2}}I(y_{k1}\leq y_{k2}+\aE{\theta}))e^{-r_{kk}+\sqrt{2r_{kk}}z_{k}},$$
which completes the proof of (\ref{eq2.1.1}).

$iii)$ We proceed as for the proof of cases $i)$ and $ii)$ using the bound \eqref{case2}.
By Lemmas 2 and 3 in \cite{Pit2004} and (\ref{eq3.4.2}), (\ref{eq3.4.3}) we obtain
\begin{eqnarray*}
A_{1}&\sim & \LE{T^{a-1}}\sum_{k=1}^{p}(e^{-x_{k}}+e^{-y_{k1}}+e^{-y_{k2}})e^{-r_{kk}+\sqrt{2r_{kk}}z_{k}}, \quad T\rightarrow\infty.
\end{eqnarray*}
With similar argument as for $A_{2}$ in the proof of case $i)$, we conclude that
$$A_{k}=o(A_{1}),\ \ \ k=2,3,4,5,6,7.$$
Further, Lemma 2 in \cite{Pit2004}  implies $A_{8}=o(A_{1})$ and
 \nelem{lem:A2} \aE{yields}
\begin{eqnarray*}
A_{10}= o(\LE{T^{a-1}})=o(A_{1}), \quad T\rightarrow\infty.
\end{eqnarray*}
\jE{Similar arguments} as for $A_{11}$ in the proof of case $ii)$ imply
$$A_{11}=o(A_{1}), \quad T\to \IF.$$
Borrowing the arguments of \aE{\cite{Pit96}},  p.\ 176  and using Lemma 3 in \cite{Pit2004} it follows that
\begin{eqnarray*}
A_{9}&=&\sum_{k=1}^{p}\mathbb{P}\left\{ \max_{t\in
[0,S]}\eta_{k}(t)> u(x_{k},z_{k}), \max_{t\in\mathfrak{R}(\delta_{2})\cap [0,S]}\eta_{k}(t)>
u(y_{k2},z_{k})\right\}\nonumber\\
&\sim& \LE{T^{a-1}}\sum_{k=1}^{p}H_{D_{2},\alpha}^{\ln H_{\alpha}+x_{k},\ln H_{D_{2},\alpha}+y_{k2}}\Quan{e^{-r_{kk}+\sqrt{2r_{kk}}z_{k}}}, \quad T\rightarrow\infty.
\end{eqnarray*}
Consequently, as $T\to \IF$
$$n_{T} P\{ \cup_{k=1}^{3p} \AK \} \sim \sum_{k=1}^{p}
(e^{-x_k} + e^{-y_{k1}}+ e^{-y_{k2}}-H_{D_{2},\alpha}^{\ln H_{\alpha}+x_{k},\ln H_{D_{2},\alpha}+y_{k2}})
e^{-r_{kk}+\sqrt{2r_{kk}}z_{k}},$$
which completes the proof of the claim in (\ref{eq2.1.1}).

$iv)$. By Lemma 5 in  \cite{Pit2004}, we have
\begin{eqnarray*}
&&\bigg|\pk{ \widetilde{\M}(T)\le \x, \widetilde{\M}(\delta_1,T)\le \y_1,\widetilde{\M}(\delta_2,T)\le \y_2}
-\pk{ \widetilde{\M}(T)\le \x, \widetilde{\M}(\delta_1,T)\le \y_1,\widetilde{\M}(T)\le \y_2}\bigg|\\
&&\leq\bigg|\pk{\widetilde{\M}(\delta_2,T)\le \y_2}
-\pk{\widetilde{\M}(T)\le \y_2}\bigg|
\rightarrow0, \quad T\rightarrow\infty.
\end{eqnarray*}
Now, by Theorem 2.1 of \cite{TanH2014}, we have
\begin{eqnarray*}
\pk{ \widetilde{\M}(T)\le \x, \widetilde{\M}(\delta_1,T)\le \y_1,\widetilde{\M}(T)\le \y_2}
&=&\pk{ \widetilde{\M}(T)\le \min(\x,\y_2), \widetilde{\M}(\delta_1,T)\le \y_1}\\
&\rightarrow&\EE{\exp\Bigl(- \sum_{k=1}^{p}
f(x_k,y_{k1},y_{k2}) e^{-r_{kk}+\sqrt{2r_{kk}}Z_{k}}\Bigr)},
\end{eqnarray*}
as $T \to \IF$ with
$$f(x_k,y_{k1},y_{k2})=
 e^{-\min(x_k,y_{k2})}+e^{-y_{k1}}$$
establishing the proof.
\hfill $\Box$

\prooftheo{Th2}
$i)$ The limiting properties of the two constants can be found in \nelem{lem:A3}. We give the proof of the relation of \eqref{eqAAA}. As for the proof of Theorem 2.1,
in view of Lemmas \Tan{3.1-3.4} and the dominated convergence theorem in order to establish the proof we need to show that	
(\ref{eq2.1.1}) holds with
\BQNY
f(x_k,y_{k1},y_{k2})&=&
e^{-x_k} + e^{-y_{k1}}+e^{-y_{k2}} - H_{D_1,\alpha}^{\ln H_\alpha+x_k,\ln H_{D_1,\alpha}+ y_{k1}} -
H_{D_2,\alpha}^{\ln H_\alpha+x_k,\ln H_{D_2,\alpha}+y_{k2} } \\
&&- H_{D_1,D_2,\alpha}^{
\ln H_{D_1,\alpha}+y_{k1},\ln H_{D_2,\alpha}+y_{k2} } +  H_{D_1,D_2,\alpha}^{\ln H_{\alpha}+x_k,\ln H_{D_1,\alpha}+y_{k1},\ln H_{D_2,\alpha}+ y_{k2}}.
\EQNY
We proceed as in the proof of case $ii)$ of Theorem 2.1 using the bound \eqref{case2}; we have thus
\BQN\label{caseT21}
\pk{ \cup_{k=1}^{3p} \AK }
&=&  \sum_{k=1}^{3p} \pk{\AK } - \sum_{1 \le k ,l \le 3p} \pk{ \AK, \AL }+ \sum_{1 \le k , l, j \le 3p} \pk{ \AK, \AL, \AJ }+\sum_{1 \le k ,\cdots, l \le 3p}\Tan{\pk{\cdot }} \notag\\
&=:& \Sigma_{1}-\Sigma_{2}+\Sigma_{3}+\Sigma_{4}.
 \EQN

By Lemmas 2 and 3 \aE{in} \cite{Pit2004} and (\ref{eq3.4.2}), (\ref{eq3.4.3}) we obtain that
\begin{eqnarray*}
\Sigma_{1}&\sim & \LE{T^{a-1}}\sum_{k=1}^{p}
(e^{-x_{k}}+e^{-y_{k1}}+e^{-y_{k2}})e^{-r_{kk}+\sqrt{2r_{kk}}z_{k}}, \quad T\rightarrow\infty.
\end{eqnarray*}
Further, write
\BQN\label{caseT22}
\Sigma_{2}=A_{2}+A_{3}+A_{4}+2A_{5}+2A_{6}+2A_{7}+A_{8}+A_{9}+A_{10},
\EQN
where $A_{i}, i=2,\cdots,10$ are defined in the proof of $ii)$ of Theorem 2.1. Hence, with similar arguments as above
$A_{i}=o(A_1)$, $i=1,\cdots,7$ and
$$A_{8}\sim \LE{T^{a-1}}\sum_{k=1}^{p}H_{D_{1},\alpha}^{\ln H_{\alpha}+x_{k},\ln H_{D_{1},\alpha}+y_{k1}}\Quan{e^{-r_{kk}+\sqrt{2r_{kk}}z_{k}}},$$
$$A_{9}\sim  \LE{T^{a-1}} \sum_{k=1}^{p}H_{D_{2},\alpha}^{\ln H_{\alpha}+x_{k},\ln H_{D_{2},\alpha}+y_{k2}}\Quan{e^{-r_{kk}+\sqrt{2r_{kk}}z_{k}}},$$
$$A_{10}\sim \LE{T^{a-1}}\sum_{k=1}^{p} H_{D_{1},D_{2}}^{\ln H_{D_{1},\alpha}+y_{k1},\ln H_{D_{2},\alpha}+y_{k2}}\Quan{e^{-r_{kk}+\sqrt{2r_{kk}}z_{k}}}$$
 as $T\rightarrow\infty$, where for the estimates of $A_{8}$ and $A_{9}$ we \jE{applied} Lemma 3 in \cite{Pit2004} and for the estimate of $A_{10}$ we have used \nelem{lem:A3}. Further 
\BQN\label{caseT23}
\Sigma_{3}& = &  \sum_{1 \le k < l<j \le 3p\atop l\neq k+p, j\neq l+p, j\neq k+2p} \mathbb{P}\{ \AK, \AL, \AJ \} + \sum_{1 \le k < l<j \le 3p\atop l= k+p, j\neq l+p, j\neq k+2p} \mathbb{P}\{ \AK, \AL, \AJ \} +\sum_{1 \le k < l<j \le 3p\atop l\neq k+p, j= l+p, j\neq k+2p} \mathbb{P}\{ \AK, \AL, \AJ \}\notag\\
 &&+ \sum_{1 \le k < l<j \le 3p\atop l\neq k+p, j\neq l+p, j= k+2p} \mathbb{P}\{ \AK, \AL, \AJ \}+\sum_{1 \le k < l<j \le 3p\atop l= k+p, j= l+p} \mathbb{P}\{ \AK, \AL, \AJ \}\notag\\
&=:&B_{1}+B_{2}+B_{3}+B_{4}+B_{5}.\notag
\EQN
For $B_{1}$, by the independence of $\eta_{k}(t)$ and $\eta_{l}(t)$,
$k\neq l$, Lemma 3 of  \cite{Pit2004} and (\ref{eq3.4.2}),
(\ref{eq3.4.3}),
 we have for some canstant \Quan{$K>0$}
\begin{eqnarray*}
B_{1}&=&\sum_{1 \le k < l<j \le 3p\atop l\neq k+p, j\neq l+p, j\neq k+2p} \mathbb{P}\{ \AK\}P\{\AL\}P\{\ \AJ \} \sim \Quan{KT^{3(a-1)}}=o(A_{1}).
\end{eqnarray*}
Similarly, we can show that
$$B_{i}\sim \Quan{KT^{2(a-1)}}=o(A_{1}),\ \ i=2,3,4.$$
For $B_{5}$, using \nelem{lem:A3}, we have
$$B_{5}\sim \LE{T^{a-1}}\sum_{k=1}^{p}H_{D_{1},D_{2},\alpha}^{\ln H_{\alpha}+x_{k},\ln H_{D_{1},\alpha}+y_{k1},\ln H_{D_{2},\alpha}+y_{k2}}e^{-r_{kk}+\sqrt{2r_{kk}}z_{k}}.$$
Finally, it is easy to see that $\Sigma_{4}=o(A_{1})$ as $T\rightarrow\infty$.
Thus, we have as $T\to \IF$
$$n_{T} \pk{ \cup_{k=1}^{3p} \AK } \sim \sum_{k=1}^{p}f(x_k,y_{k1},y_{k2})e^{-r_{kk}+\sqrt{2r_{kk}}z_{k}},$$
with
\BQNY
f(x_k,y_{k1},y_{k2})&=&
e^{-x_k} + e^{-y_{k1}}+e^{-y_{k2}} - H_{D_1,\alpha}^{\ln H_\alpha+x_k,\ln H_{D_1,\alpha}+ y_{k1}} -
H_{D_2,\alpha}^{\ln H_\alpha+x_k,\ln H_{D_2,\alpha}+y_{k2} } \\
&&- H_{D_1,D_2,\alpha}^{
\ln H_{D_1,\alpha}+y_{k1},\ln H_{D_2,\alpha}+y_{k2} } +  H_{D_1,D_2,\alpha}^{\ln H_{\alpha}+x_k,\ln H_{D_1,\alpha}+y_{k1},\ln H_{D_2,\alpha}+ y_{k2}},
\EQNY
which completes the proof of (\ref{eq2.1.1}).

$ii)$ Applying Lemma 5 in \cite{Pit2004} we obtain 
\begin{eqnarray*}
&&\bigg|\pk{ \widetilde{\M}(T)\le \x, \widetilde{\M}(\delta_1,T)\le \y_1,\widetilde{\M}(\delta_2,T)\le \y_2}
-\pk{ \widetilde{\M}(T)\le \x, \widetilde{\M}(\delta_1,T)\le \y_1,\widetilde{\M}(T)\le \y_2}\bigg|\\
&&\leq\bigg|\pk{\widetilde{\M}(\delta_2,T)\le \y_2}
-\pk{\widetilde{\M}(T)\le \y_2}\bigg|
\rightarrow0, \quad T\rightarrow\infty.
\end{eqnarray*}
Further, Theorem 2.2 in \cite{TanH2014} yields 
\begin{eqnarray*}
\pk{ \widetilde{\M}(T)\le \x, \widetilde{\M}(\delta_1,T)\le \y_1,\widetilde{\M}(T)\le \y_2}
&=&\pk{ \widetilde{\M}(T)\le \min(\x,\y_2), \widetilde{\M}(\delta_1,T)\le \y_1}\\
&\rightarrow&\EE{\exp\Bigl(- \sum_{k=1}^{p}
f(x_k,y_{k1},y_{k2}) e^{-r_{kk}+\sqrt{2r_{kk}}Z_{k}}\Bigr)}
\end{eqnarray*}
with
$$f(x_k,y_{k1},y_{k2})=
e^{-\min(x_k,y_{k2})}+ e^{-y_{k1}} - H_{D_1,\alpha}^{\ln H_\alpha+\min(x_k,y_{k2}),\ln H_{D_1,\alpha}+ y_{k1}},$$
which completes the proof.

$iii)$ By Theorem 2.3 in \cite{TanH2014} for the dense grid $\RDI  ,i=1,2$ and any $\y_i\inr^p$
\BQNY
\limit{T} \pk{\widetilde{\M}(\iH{\delta_i},T) \le \y_i} = \E{\exp\Bigl( - \sum_{k=1}^p e^{- y_{ki} - r_{kk} + \sqrt{2 r_{kk}} Z_k} \Bigr)}
\EQNY
and further
\BQNY
\limit{T} \pk{\widetilde{\M}(T) \le \x} = \E{\exp\Bigl( - \sum_{k=1}^p e^{- \Quan{x_{k}} - r_{kk} + \sqrt{2 r_{kk}} Z_k} \Bigr)},
\quad \forall \x\inr^p,
\EQNY
 hence  the claim follows immediately from \nelem{lem:A4}.\QED

\section{Appendix}
For the proof of the main results, we need the following technical lemmas. Let in the sequel
$\mathcal{C}$ be \bE{a} \jE{positive} constant whose value will change from place to place and $\overline{\Phi}, \varphi$ be the survival function and the density function of \bE{an $N(0,1)$ random} variable, respectively.

\BL\label{lem:A1} Suppose that
$\mathfrak{R}(\delta_{1})$ and $\mathfrak{R}(\delta_{2})$ are  sparse grids \AH{and $a\in (0,1)$}.\\
i) If $\limit{T} \delta_1(T)/ \delta_2(T)=\infty$ or $\mathfrak{R(\delta_{1})}\cap\mathfrak{R(\delta_{2})}=\emptyset$, then we have \Tan{for $k\leq p$} as $T\rightarrow\infty$
\begin{eqnarray*}
&&\mathbb{P}\left\{\max_{t\in \mathcal{R}(\delta_{1})\cap[0,S]}\eta_{k}(t)> u(y_{k1},z_{k}), \max_{t\in \mathcal{R}(\delta_{2})\cap[0,S]}\eta_{k}(t)> u(y_{k2},z_{k})\right\}=o(\LE{T^{a-1}}).
\end{eqnarray*}
ii) Let $\mathfrak{R(\delta_{1})}\cap\mathfrak{R(\delta_{2})}=\mathfrak{R(\delta_{3})}$ and
$\lim_{T\to \IF} \ln(\frac{\delta_{3}(T)}{\delta_{1}(T)})= \theta_{1} \in [0,\IF)$, $\lim_{T\to \IF} \ln(\frac{\delta_{3}(T)}{\delta_{2}(T)})= \theta_{2} \in [0,\IF)$ hold. If $y_{k1}> y_{k2}+ \theta_{2}-\theta_{1}$, then we have \Tan{for $k\leq p$} as $T\rightarrow\infty$
\begin{eqnarray*}
&&\mathbb{P}\left\{\max_{t\in \mathcal{R}(\delta_{1})\cap[0,S]\setminus \mathcal{R}(\delta_{3})\cap[0,S]}\eta_{k}(t)> u(y_{k1},z_{k}), \max_{t\in \mathcal{R}(\delta_{2})\cap[0,S]}\eta_{k}(t)> u(y_{k2},z_{k})\right\}=o(\LE{T^{a-1}}),
\end{eqnarray*}
whereas \iE{if} $y_{k1}\leq y_{k2}+ \theta_{2}-\theta_{1}$ 
\begin{eqnarray*}
&&\mathbb{P}\left\{\max_{t\in \mathcal{R}(\delta_{1})\cap[0,S]}\eta_{k}(t)> u(y_{k1},z_{k}), \max_{t\in \mathcal{R}(\delta_{2})\cap[0,S]\setminus \mathcal{R}(\delta_{3})\cap[0,S]}\eta_{k}(t)> u(y_{k2},z_{k})\right\}=o(\LE{T^{a-1}})
\end{eqnarray*}
holds.
\EL

\prooflem{lem:A1} The following fact will be extensively used in the proof. From assumption (\ref{eq1.3}), we can choose \iE{an} $\epsilon>0$ such that for all
$|s-t|\leq\epsilon<2^{-1/\alpha}$
\begin{eqnarray}
\label{eq.A14}
\frac{1}{2}|s-t|^{\alpha}\leq 1-r_{kk}(s,t)\leq 2|s-t|^{\alpha}.
\end{eqnarray}
$i)$ We first deal with the case $\lim_{T\rightarrow\infty}\delta_{1}(T)/\delta_{2}(T)=\infty$.
It is easy to check that
\begin{eqnarray*}
&&\sum_{k=1}^{p}\mathbb{P}\left\{\max_{t\in \mathcal{R}(\delta_{1})\cap[0,S]}\eta_{k}(t)> u(y_{k1},z_{k}), \max_{t\in \mathcal{R}(\delta_{2})\cap[0,S]}\eta_{k}(t)> u(y_{k2},z_{k})\right\}\nonumber\\
&\iH{\leq} &\sum_{k=1}^{p}\bigg[\mathbb{P}\left\{\max_{t\in \mathcal{R}(\Quan{\delta_{1}})\cap[0,S]}\eta_{k}(t)> u(y_{k2},z_{k})\right\}\\
&&+ \mathbb{P}\left\{\max_{t\in \mathcal{R}(\delta_{1})\cap[0,S]}\eta_{k}(t)> u(y_{k1},z_{k}), \max_{t\in \mathcal{R}(\delta_{2})\cap[0,S]\setminus \mathcal{R}(\delta_{1})\cap[0,S]}\eta_{k}(t)> u(y_{k2},z_{k})\right\}\bigg].
\end{eqnarray*}
By Lemma 2 of \cite{Pit2004} and the definition of $ u(y_{k2},z_{k})$,
 we have as $T\rightarrow\infty$
\BQNY
\mathbb{P}\left\{\max_{t\in\mathcal{R}(\delta_{1})\cap[0,S]}\eta_{k}(t)> u(y_{k2},z_{k})\right\}
&\sim &S\delta_{1}^{-1}(T)\overline{\Phi}( u(y_{k2},z_{k}))\\
&=& \mathcal{C}S\delta_{1}^{-1}(T)T^{-1}\delta_{2}(T)\\
&=&\mathcal{C}T^{a-1}\frac{\delta_{2}(T)}{\delta_{1}(T)}=o(T^{a-1}).
\EQNY
Now, for $m, n\in \mathbb{N}$ and the $\epsilon$ chosen in (\ref{eq.A14}), we have
\BQN
\label{eq.A11}
&&\mathbb{P}\left\{\max_{t\in \mathcal{R}(\delta_{1})\cap[0,S]}\eta_{k}(t)> u(y_{k1},z_{k}), \max_{t\in \mathcal{R}(\delta_{2})\cap[0,S]}\eta_{k}(t)> u(y_{k2},z_{k})\right\}\nonumber\\
&&= \mathbb{P}\left\{\max_{t\in \mathcal{R}(\delta_{1})\cap[0,S]}\eta_{k}(t)> u(y_{k1},z_{k}), \max_{t\in \mathcal{R}(\delta_{2})\cap[0,S]\setminus \mathcal{R}(\delta_{1})\cap[0,S]}\eta_{k}(t)> u(y_{k2},z_{k})\right\}+o(T^{a-1})\nonumber\\
&&\Quan{\leq\sum_{n=0}^{[S/\delta_{1}]+1}\mathbb{P}\left\{\eta_{k}(n\delta_{1})> u(y_{k1},z_{k}), \max_{t\in \mathcal{R}(\delta_{2})\cap[0,S]\setminus \mathcal{R}(\delta_{1})\cap[0,S]}\eta_{k}(t)>  u(y_{k2},z_{k})\right\}+o(T^{a-1})}\nonumber\\
&&\Quan{\leq \sum_{n=0}^{[S/\delta_{1}]+1}\mathbb{P}\left\{\eta_{k}(n\delta_{1})> u(y_{k1},z_{k}), \max_{0\leq t\leq S \atop |t-n\delta_{1}|\leq\epsilon}\eta_{k}(t)> u(y_{k2},z_{k})\right\}}\nonumber\\
&&\ \ \ \ \Quan{ + \sum_{n=0}^{[S/\delta_{1}]+1}\mathbb{P}\left\{\eta_{k}(n\delta_{1})> u(y_{k1},z_{k}), \max_{0\leq m\delta_{2}\leq S \atop |n\delta_{1}-m\delta_{2}|>\epsilon}\eta_{k}(m\delta_{2})> u(y_{k2},z_{k})\right\}+o(T^{a-1})}\nonumber\\
&&=: S_{T,1}+S_{T,2}+o(T^{a-1}),\nonumber
\EQN
\Quan{where $[x]$ denotes the integer part of $x$.
\jE{By stationarity we have setting  $\eta_{nk}^*(t)= \eta_k(n\delta_{1})+ \eta_k(t)$}
\begin{eqnarray*}
S_{T,1}&\leq&  \sum_{n=0}^{[S/\delta_{1}]+1}\mathbb{P}\left\{\max_{n\delta_{1}-\epsilon\leq t\leq n\delta_{1}+ \epsilon}\jE{\eta_{nk}^*}(t)> u(y_{k1},z_{k})+ u(y_{k2},z_{k})\right\}\\
&=&\mathcal{C} \frac{S}{\delta_{1}}\mathbb{P}\left\{\max_{0\leq t< \epsilon}\jE{\eta^*_{0k}}(t)> u(y_{k1},z_{k})+ u(y_{k2},z_{k})\right\}.
\end{eqnarray*}}
For the correlation function of $\jE{\eta^*_{0k}(t)}=\eta_{k}(0)+\eta_{k}(t)$, $t\in [0,\epsilon]$ we have
\begin{eqnarray*}
1-\frac{\mathbb{E}(\jE{\eta^*_{0k}(s)})(\jE{\eta^*_{0k}(t)})}
{\sqrt{    \mathbb{E}((\jE{\eta^*_{0k}(s)})^{2})
\mathbb{E}((\jE{\eta^*_{0k}(t)})^{2} )} }
& \leq& \frac{1-r_{kk}(t-s)}{2\sqrt{1+r_{kk}(t)}\sqrt{1+r_{kk}(s)}}\\
&\leq& \frac{2|t-s|^{\alpha}}{2-2\epsilon^{\alpha}}\leq 1-\exp(-|t-s|^{\alpha}).
\end{eqnarray*}
Further 
$$Var(\jE{\eta^*_{0k}(t)})=2+2r_{kk}(t)=4-2|t|^{\alpha}(1+o(1))$$
as $t\rightarrow 0$. Hence by Slepian's inequality (see e.g. Theorem 7.4.2 of \cite{leadbetter1983extremes}) we have
\begin{eqnarray*}
&&\mathbb{P}\left\{\max_{0\leq t< \epsilon}\jE{\eta^*_{0k}(t)}> u(y_{k1},z_{k})+ u(y_{k2},z_{k})\right\}\\
&&= \mathbb{P}\left\{\max_{0\leq t< \epsilon}\frac{\jE{\eta^*_{0k}(t)}}
{\sqrt{    \mathbb{E}((\jE{\eta^*_{0k}(t)})^{2})}}
\sqrt{    \mathbb{E}((\jE{\eta^*_{0k}(t)})^{2})}> u(y_{k1},z_{k})+ u(y_{k2},z_{k})\right\}\\
&& \leq  \mathbb{P}\left\{\max_{0\leq t< \epsilon}\aE{W(t)}\sqrt{    \mathbb{E}((\jE{\eta^*_{0k}(t)})^{2})}> u(y_{k1},z_{k})+ u(y_{k2},z_{k})\right\},
\end{eqnarray*}
where $\aE{W}$ is a Gaussian zero mean stationary process with covariance function $\exp(-|t|^{\alpha})$, thus the condition of Theorem D.3 \aE{
in \cite{Pit96}} for the case $\alpha=\beta$ hold. By \aE{that} theorem
\begin{eqnarray*}
S_{T,1}&\leq& \mathcal{C }\frac{S}{\delta_{1}} \overline{\Phi}\left(\frac{u(y_{k1},z_{k})+ u(y_{k2},z_{k})}{2}\right).
\end{eqnarray*}
The definition of $u(y_{ki},z_{k}), i=1,2$ implies thus for sparse grids
\begin{eqnarray}
\label{eq.A13}
&&[u(y_{ki},z_{k})]^{2}=2\ln T-\ln\ln T+2\ln \delta_{i}^{-1}(T)+O(1).
\end{eqnarray}
Consequently, from the fact that $\limit{T} \delta_1(T)/ \delta_2(T)=\infty$
$$S_{T,1}\leq \mathcal{ C }\frac{S}{\delta_{1}(T)}\frac{1}{\sqrt{\ln T}}T^{-1}\sqrt{\ln T}\delta_{1}^{1/2}(T)\delta_{2}^{1/2}(T)
= \mathcal{ C }T^{a-1}\left(\frac{\delta_{2}(T)}{\delta_{1}(T)}\right)^{1/2}=o(T^{a-1}), \quad T \to \IF.$$
Now, let $\vartheta_{kk}(t)=\sup_{t\leq s\leq S}r_{kk}(s)$.
Assumption (\ref{eq1.3}) implies that
$\vartheta_{kk}(\epsilon)<1$ for all $T$ and any
$\epsilon\in(0,2^{-1/\alpha})$. Consequently, we may choose some
positive constant $\beta_{kk}$ such that
$$
\beta_{kk}<\frac{1-\vartheta_{kk}(\epsilon)}{1+\vartheta_{kk}(\epsilon)} < 1
$$
for all sufficiently large $T$. In the following we choose
\BQN \label{bkk}
0 < a <b<\min_{1 \le k \le p }\beta_{kk}.
 \EQN
For the second term, by stationarity and Berman's inequality (see eg. Theorem 4.2.1 of \cite{leadbetter1983extremes}, Theorem C.2 of \cite{Pit96}), we have
\begin{eqnarray*}
S_{T,2}&\leq& \sum_{n=0}^{[S/\delta_{1}]+1}\sum_{0\leq m\delta_{2}\leq S\atop |n\delta_{1}-m\delta_{2}|>\epsilon}\mathbb{P}\left\{\eta_{k}(n\delta_{1})> u(y_{k1},z_{k}), \eta_{k}(m\delta_{2})> u(y_{k2},z_{k})\right\}\\
&\leq & \sum_{n=0}^{[S/\delta_{1}]+1}\sum_{0\leq m\delta_{2}\leq S\atop |n\delta_{1}-m\delta_{2}|>\epsilon}\bigg[\overline{\Phi}( u(y_{k1},z_{k}))
\overline{\Phi}( u(y_{k2},z_{k}))+\aE{\mathcal{C}}\exp\left(\frac{u^{2}(y_{k1},z_{k})+ u^{2}(y_{k2},z_{k})}{2(1+r_{kk}(|n\delta_{1}-m\delta_{2}|))}\right)\bigg]\\
&\leq & \frac{S}{\delta_{1}}\frac{S}{\delta_{2}}\bigg[\overline{\Phi}( u(y_{k1},z_{k}))\overline{\Phi}( u(y_{k2},z_{k}))
+\aE{\mathcal{C}}\exp\left(\frac{u^{2}(y_{k1},z_{k})+ u^{2}(y_{k2},z_{k})}{2(1+\vartheta_{kk}(\epsilon))}\right)\bigg]\\
&=:&S_{T,21}+S_{T,22}.
\end{eqnarray*}
Utilising again (\ref{eq.A13}) 
\begin{eqnarray*}
S_{T,21}&\leq&  \mathcal{C}\frac{S}{\delta_{1}}\frac{S}{\delta_{2}}\frac{\varphi( u(y_{k1},z_{k}))}{ u(y_{k1},z_{k})}\frac{\varphi( u(y_{k2},z_{k}))}{ u(y_{k2},z_{k})}\\
&\leq&  \mathcal{C}\frac{S}{\delta_{1}}\frac{S}{\delta_{2}}\frac{1}{\ln T}\exp\left(-\frac{1}{2}u^{2}(y_{k1},z_{k})\right)\exp\left(-\frac{1}{2}u^{2}(y_{k2},z_{k})\right)\\
&\leq& \mathcal{C}\frac{S}{\delta_{1}}\frac{S}{\delta_{2}}\frac{1}{\ln T}T^{-1}(\ln T)^{1/2}\delta_{1}T^{-1}(\ln T)^{1/2}\delta_{2}\\
&=&\mathcal{C}T^{2(a-1)}
\end{eqnarray*}
as $T\rightarrow\infty$.
Since  $u(y_{ki},z_{k})\sim (2\ln T)^{1/2},i=1,2$
\begin{eqnarray*}
\label{eq.A16}
S_{T,22}&\leq&\mathcal{C}\frac{S}{\delta_{1}}\frac{S}{\delta_{2}}\exp\left(\frac{u^{2}(y_{k1},z_{k})+ u^{2}(y_{k2},z_{k})}{2(1+\vartheta_{kk}(\epsilon))}\right)\\
&\leq & \mathcal{C} \frac{T^{a}}{\delta_{1}}\frac{T^{a}}{\delta_{2}}T^{-\frac{2}{1+\vartheta_{kk}(\epsilon)}}\nonumber\\
&\leq & \mathcal{C} T^{a-1} T^{a-\frac{1-\vartheta_{_{kk}}(\epsilon)}{1+\vartheta_{kk}(\epsilon)}}(\delta_{1}\delta_{2})^{-1}.
\end{eqnarray*}
Both (\ref{bkk}) and  $\limit{T} (\ln T)^{1/\alpha}\delta_{i}(T)=\infty$ imply
$ S_{T,22}=o(T^{a-1})$ as $T\rightarrow\infty$.\\
Let us consider now the case that $\mathfrak{R(\delta_{1})}\cap\mathfrak{R(\delta_{2})}=\emptyset$.
Without loss of generality, we suppose that $u(y_{k1},z_{k})< u(y_{k2},z_{k})$ holds for sufficient large $T$.
By stationarity, for $m, n\in \mathbb{N}$ and  $\epsilon>0$ we have
\begin{eqnarray*}
\label{eq.A11}
&&\mathbb{P}\left\{\max_{t\in \mathcal{R}(\delta_{1})\cap[0,S]}\eta_{k}(t)> u(y_{k1},z_{k}), \max_{t\in \mathcal{R}(\delta_{2})\cap[0,S]}\eta_{k}(t)> u(y_{k2},z_{k})\right\}\nonumber\\
&&\leq\sum_{n=0}^{[S/\delta_{1}]+1}\mathbb{P}\left\{\eta_{k}(n\delta_{1})> u(y_{k1},z_{k}), \max_{0< m\delta_{2}\leq S}\eta_{k}(m\delta_{2})> u(y_{k2},z_{k})\right\}\nonumber\\
&&\leq \sum_{n=0}^{[S/\delta_{1}]+1}\mathbb{P}\left\{\eta_{k}(n\delta_{1})> u(y_{k1},z_{k}), \max_{0< m\delta_{2}\leq S\atop |n\delta_{1}-m\delta_{2}|\leq \epsilon}\eta_{k}(m\delta_{2})> u(y_{k1},z_{k})\right\}\nonumber\\
&&\ \ \ \ + \sum_{n=0}^{[S/\delta_{1}]+1}\mathbb{P}\left\{\eta_{k}(n\delta_{1})> u(y_{k1},z_{k}), \max_{0< m\delta_{2}\leq S\atop |n\delta_{1}-m\delta_{2}|> \epsilon}\eta_{k}(m\delta_{2})> u(y_{k2},z_{k})\right\}\nonumber\\
&&= \mathcal{C} \frac{S}{\delta_{1}}\mathbb{P}\left\{\eta_{k}(0)> u(y_{k1},z_{k}), \max_{0< m\delta_{2}\leq \epsilon}\eta_{k}(m\delta_{2})> u(y_{k1},z_{k})\right\}\nonumber\\
&&\ \ \ \ + \sum_{n=0}^{[S/\delta_{1}]+1}\mathbb{P}\left\{\eta_{k}(n\delta_{1})> u(y_{k1},z_{k}), \max_{0< m\delta_{2}\leq S\atop |n\delta_{1}-m\delta_{2}|> \epsilon}\eta_{k}(m\delta_{2})> u(y_{k2},z_{k})\right\}\nonumber\\
&&=: R_{T,1}+R_{T,2}.
\end{eqnarray*}
\iH{Using \Quan{the} well-known results for bivariate Gaussian tail probability (see e.g., \cite{EH05}) setting} $r=r_{kk}(m\delta_{2})$ we have 
\begin{eqnarray*}
\label{eq.A12}
R_{T,1}&\iH{\le} &\mathcal{C} \frac{S}{\delta_{1}}\sum_{0< m\delta_{2}\leq \epsilon}\mathbb{P}\left\{\eta_{k}(0)> u(y_{k1},z_{k}), \eta_{k}(m\delta_{2})> u(y_{k1},z_{k})\right\}\nonumber\\
\COM{&=& \sum_{n=0}^{S/\delta_{1}}\sum_{0<m\delta_{2}\leq \epsilon}\int_{u(y_{k1},z_{k})}^{\infty}\mathbb{P}\left\{\eta_{k}(m\delta_{2})> u(y_{k1},z_{k})|\eta_{k}(0)=x\right\}\iH{\varphi}(x)dx\nonumber\\
&=& \frac{S}{\delta_{1}}\sum_{0<m\delta_{2}\leq \epsilon}\int_{u(y_{k1},z_{k})}^{\infty}
\overline{\Phi}\left(\frac{u(y_{k1},z_{k})-xr}{\sqrt{1-r^{2}}}\right) d \varphi(x) dx\\
&=&\frac{S}{\delta_{1}}\sum_{0<m\delta_{2}\leq \epsilon}\left[-\overline{\Phi}(x)\overline{\Phi}\left(\frac{u(y_{k1},z_{k})-xr}{\sqrt{1-r^{2}}}\right)\right]_{x=u(y_{k1},z_{k})}^{\infty}\\
&&+\frac{S}{\delta_{1}}\sum_{0<m\delta_{2}\leq \epsilon}\int_{u(y_{k1},z_{k})}^{\infty}
 \overline{\Phi}(u(y_{k1},z_{k}))\frac{r}{\sqrt{1-r^{2}}}\varphi\left(\frac{u(y_{k1},z_{k})-xr}{\sqrt{1-r^{2}}}\right)dx\\
&\leq& \frac{S}{\delta_{1}}\sum_{0<m\delta_{2}\leq \epsilon}\left[\overline{\Phi}(u(y_{k1},z_{k}))
\overline{\Phi}\left(\frac{u(y_{k1},z_{k})-u(y_{k1},z_{k})r}{\sqrt{1-r^{2}}}\right)\right]\\
&&+\frac{S}{\delta_{1}}\sum_{0<m\delta_{2}\leq \epsilon}\int_{u(y_{k1},z_{k})}^{\infty}\overline{\Phi}(u(y_{k1},z_{k}))
\frac{r}{\sqrt{1-r^{2}}}\varphi\left(\frac{u(y_{k1},z_{k})-xr}{\sqrt{1-r^{2}}}\right)dx\\
}
&=& \iH{\mathcal{C} \frac{S}{\delta_{1}}\sum_{0< m\delta_{2}\leq \epsilon}\left[\overline{\Phi}(u(y_{k1},z_{k}))
\overline{\Phi}\left(u(y_{k1},z_{k})\frac{\sqrt{1-r}}{\sqrt{1+r}}\right)\right].}
\end{eqnarray*}
Since by (\ref{eq.A14})
$$\frac{1-r}{1+r}=\frac{1-r_{kk}(m\delta_{2})}{1+r_{kk}(m\delta_{2})}\geq \frac{1}{4}(m\delta_{2})^{\alpha}$$
and using  (\ref{eq.A13}) we obtain
\begin{eqnarray*}
\label{eq.A12}
R_{T,1}&\leq&\mathcal{C}  \frac{S}{\delta_{1}}\sum_{0< m\delta_{2}\leq \epsilon}\left[\overline{\Phi}(u(y_{k1},z_{k}))
\overline{\Phi}\left(\frac{1}{2}(m\delta_{2})^{\alpha/2}u(y_{k1},z_{k})\right)\right]\\
&=& \mathcal{C}T^{a-1}\sum_{0< m\delta_{2}\leq \epsilon}\overline{\Phi}\left(\frac{1}{2}(m\delta_{2})^{\alpha/2}u(y_{k1},z_{k})\right)\\
&=& \mathcal{C} T^{a-1}\sum_{0< m\delta_{2}\leq \epsilon}\frac{1}{(m\delta_{2})^{\alpha/2}u(y_{k1},z_{k})}\exp\left(-\frac{1}{8}(m\delta_{2})^{\alpha}u^{2}(y_{k1},z_{k})\right)\\
&=&\mathcal{C} T^{a-1}\sum_{0< m\delta_{2}\leq \epsilon}\frac{1}{[m\delta_{2}(\ln T)^{1/\alpha}]^{\alpha/2}}\exp\left(-\frac{1}{4}[m\delta_{2}(\ln T)^{1/\alpha}]^{\alpha}\right)\\
&\leq &\mathcal{C} T^{a-1}\frac{1}{[\delta_{2}(\ln T)^{1/\alpha}]^{\alpha/2}}\sum_{0<m\leq [\epsilon/\delta_{2}]+1}\exp\left(-\frac{1}{4}[m\delta_{2}(\ln T)^{1/\alpha}]^{\alpha}\right)\\
&\leq &\mathcal{C} T^{a-1}\frac{1}{[\delta_{2}(\ln T)^{1/\alpha}]^{\alpha/2}}\\
&=&T^{a-1}o(1),
\end{eqnarray*}
\aE{where we used additionally}  the fact that $\lim_{T\rightarrow\infty}(\ln T)^{1/\alpha}\delta_{i}(T)=\infty$, $i=1,2$.
\jE{By repeating} the calculations for $S_{T,2}$ we obtain further
$R_{T,2}=o(T^{a-1})$ as $T\rightarrow\infty$, which completes the proof.

$ii)$ If $y_{k1}\leq y_{k2}+\theta_{2}-\theta_{1}$, then
we have $u(y_{k1},z_{k})\leq u(y_{k2},z_{k})$ for sufficient large $T$. \Quan{By stationarity we have} for $m, n\in \mathbb{N}$ and $\epsilon>0$
\BQN
\label{eq.A11}
&&\mathbb{P}\left\{\max_{t\in \mathcal{R}(\delta_{1})\cap[0,S]}\eta_{k}(t)> u(y_{k1},z_{k}), \max_{t\in \mathcal{R}(\delta_{2})\cap[0,S]\setminus\mathcal{R}(\delta_{3})\cap[0,S] }\eta_{k}(t)> u(y_{k2},z_{k})\right\}\nonumber\\
&&\leq \sum_{n=0}^{[S/\delta_{1}]+1}\mathbb{P}\left\{\eta_{k}(n\delta_{1})> u(y_{k1},z_{k}), \max_{t\in \mathcal{R}(\delta_{2})\cap[0,S]\setminus\mathcal{R}(\delta_{3})\cap[0,S]}\eta_{k}(t)> u(y_{k2},z_{k})\right\}\nonumber\\
&&\leq \sum_{n=0}^{[S/\delta_{1}]+1}\mathbb{P}\left\{\eta_{k}(n\delta_{1})> u(y_{k1},z_{k}), \max_{0< m\delta_{2}\leq S\atop |n\delta_{1}-m\delta_{2}|\leq\epsilon}\eta_{k}(m\delta_{2})> u(y_{k2},z_{k})\right\}\nonumber\\
&&\ \ \ + \sum_{n=0}^{[S/\delta_{1}]+1}\mathbb{P}\left\{\eta_{k}(n\delta_{1})> u(y_{k1},z_{k}), \max_{0< m\delta_{2}\leq S\atop |n\delta_{1}-m\delta_{2}|>\epsilon}\eta_{k}(m\delta_{2})> u(y_{k2},z_{k})\right\}\nonumber\\
&&\leq\mathcal{C}\frac{S}{\delta_{1}}\mathbb{P}\left\{\eta_{k}(0)> u(y_{k1},z_{k}), \max_{0< m\delta_{2}\leq \epsilon}\eta_{k}(m\delta_{2})> u(y_{k1},z_{k})\right\}\nonumber\\
&&\ \ \ \ + \sum_{n=0}^{[S/\delta_{1}]+1}\mathbb{P}\left\{\eta_{k}(n\delta_{1})> u(y_{k1},z_{k}), \max_{0< m\delta_{2}\leq S\atop |n\delta_{1}-m\delta_{2}|>\epsilon}\eta_{k}(m\delta_{2})> u(y_{k2},z_{k})\right\}\nonumber\\
&&=: M_{T,1}+M_{T,2}.\nonumber
\EQN
Using the same estimates for $R_{T,1}$ and $R_{T,2}$, we get that both $M_{T,1}$ and $M_{T,2}$ are $o(T^{a-1})$. The proof \iE{when} $y_{k1}> y_{k2}+(\theta_{2}-\theta_{1})$ is similar.
This completes the proof of the lemma.
\QED

\COM{
By Normal Comparison Lemma (see eg. Theorem 4.2.1 of \cite{leadbetter1983extremes}, Theorem C.2 of \cite{Pit96}), for $n,m\in \mathbb{N}$, it is easy to check that  both the left hand side of the two assertions is bounded above by
\begin{eqnarray}
\label{eq.A11}
\sum_{n\delta_{1}\in \mathcal{S}_{1},m\delta_{2}\in \mathcal{S}_{1}, n\delta_{1}\neq m\delta_{2}}|r_{kk}(n\delta_{1},m\delta_{2})|\exp\bigg(-\frac{(u(y_{k1},z_{k}))^{2}+(u(y_{k2},z_{k}))^{2}}{2(1+r_{kk}(n\delta_{1},m\delta_{2}))}\bigg)
\end{eqnarray}
as $T\rightarrow\infty$, where $\mathcal{S}_{1}=[0,S]=[0,T^{a}]$.

Split the sum into
two parts as for some $\epsilon>0$
\begin{eqnarray}
\label{eq.A12}
\sum_{n\delta_{1}, m\delta_{2}\in \mathcal{S}_{1}\atop  |n\delta_{1}-m\delta_{2}|<\epsilon}+\sum_{n\delta_{1}, m\delta_{1}\in \mathcal{S}_{1}\atop |n\delta_{1}-m\delta_{2}|\geq\epsilon}=:S_{T,1}+S_{T,2}.
\end{eqnarray}
The definition of $u(y_{ki},z_{k}), i=1,2$ implies thus
\begin{eqnarray}
\label{eq.A13}
u(y_{ki},z_{k})^{2}=2\ln T-\ln\ln T+2\ln \delta_{i}^{-1}+O(1), \ \ i=1,2
\end{eqnarray}
for sparse grids.
Note that Assumption (\ref{eq1.3}) implies for all
$|s-t|\leq\epsilon<2^{-1/\alpha}$
\begin{eqnarray}
\label{eq.A14}
\frac{1}{2}|s-t|^{\alpha}\leq 1-r_{kk}(s,t)\leq 2|s-t|^{\alpha}.
\end{eqnarray}
Consequently, we have (write $u_{12}=(u(y_{k1},z_{k}))^{2}+(u(y_{k2},z_{k}))^{2}$)
\begin{eqnarray*}
\label{eq.A15}
S_{T,1}
&\leq &\mathcal{C} \sum_{n\delta_{1}, m\delta_{2}\in \mathcal{S}_{1}\atop  |n\delta_{1}-m\delta_{2}|<\epsilon}
\exp\left(-\frac{u_{12}}{4}\right)\exp\left(-\frac{(1-r_{kk}(n\delta_{1},m\delta_{2}))u_{12}}{4(1+r_{kk}(n\delta_{1},m\delta_{2}))}\right)\nonumber\\
&\leq & \mathcal{C} T^{-1}(\ln T)^{1/2}\LE{(\delta_{1}\delta_{2})^{1/2}}\sum_{n\delta_{1}, m\delta_{2}\in \mathcal{S}_{1}\atop  |n\delta_{1}-m\delta_{2}|<\epsilon}\exp\left(-\frac{(1-r_{kk}(n\delta_{1},m\delta_{2}))u_{12}}{4(1+r_{kk}(n\delta_{1},m\delta_{2}))}\right)\nonumber\\
&\leq & \mathcal{C} T^{-1}(\ln T)^{1/2}\LE{(\delta_{1}\delta_{2})^{1/2}} \sum_{n\delta_{1}, m\delta_{2}\in \mathcal{S}_{1}\atop  |n\delta_{1}-m\delta_{2}|<\epsilon}\exp\left(-\frac{1}{16}|n \delta_{1}-m\delta_{2}|^{\alpha}u_{12}\right)\nonumber\\
&\leq & \mathcal{C} T^{-1}(\ln T)^{1/2}\LE{(\delta_{1}\delta_{2})^{1/2}}\sum_{n\delta_{1}, m\delta_{2}\in \mathcal{S}_{1}\atop  |n\delta_{1}-m\delta_{2}|<\epsilon}\exp\left(-\frac{1}{4}|n\delta_{1}-m\delta_{2}|^{\alpha}\ln T\right)\nonumber\\
&\leq & \mathcal{C} T^{a-1}(\ln T)^{1/2}\delta_{1}^{1/2}\delta_{2}^{-1/2}
\sum_{0<n\delta_{1}<\epsilon}\exp\left(-\frac{1}{4}(n\delta_{1})^{\alpha}\ln T\right).
\end{eqnarray*}
Since $\mathfrak{R}(\delta_{i})$ is a sparse grid we have
 $ \limit{T}\delta_{i}(\ln
T)^{1/\alpha}=\infty$. A straightforward calculation shows that
\begin{eqnarray*}
\sum_{0<n\delta_{1}<\epsilon}
\exp\left(-\frac{1}{4}(n\delta_{1})^{\alpha}\ln T\right)\leq \mathcal{C}.
\end{eqnarray*}
Hence  $S_{T,1}\leq \mathcal{C} T^{a-1}(\ln T)^{1/2}\delta_{1}^{1/2}\delta_{2}^{-1/2}=o(T^{a-1})$ by
using the fact that $\delta_{i}(\ln T)^{1/\alpha}\rightarrow\infty$.

Now, let $\vartheta_{kk}(t)=\sup_{t<|n\delta_{1}-m\delta_{2}|\leq T}\{r_{kk}(n\delta_{1},m\delta_{2})\}$.
Assumption (\ref{eq1.3}) implies that
$\vartheta_{kk}(\epsilon)<1$ for all $T$ and any
$\epsilon\in(0,2^{-1/\alpha})$. Consequently, we may choose some
positive constant $\beta_{kk}$ such that
$$
\beta_{kk}<\frac{1-\vartheta_{kk}(\epsilon)}{1+\vartheta_{kk}(\epsilon)}
$$
for all sufficiently large $T$. In the following we choose
\BQN \label{bkk}
0 < a <b<\min_{k\in\{1,\cdots,p\}}\beta_{kk}.
 \EQN
Since  $u(y_{ki})\sim (2\ln T)^{1/2},i=1,2$
\begin{eqnarray*}
\label{eq.A16}
S_{T,2}&\leq&\mathcal{C}\sum_{n\delta_{1}, m\delta_{1}\in \mathcal{S}_{1}\atop |n\delta_{1}-m\delta_{2}|\geq\epsilon}|r_{kk}(n\delta_{1},m\delta_{2})|\exp\left(-\frac{u_{12}}{2(1+\vartheta_{kk}(n\delta_{1},m\delta_{2})}\right)\nonumber\\
&\leq & \mathcal{C} \exp\left(-\frac{u_{12}}{2(1+\vartheta_{kk}(\epsilon))}\right)\sum_{n\delta_{1}, m\delta_{1}\in \mathcal{S}_{1}\atop |n\delta_{1}-m\delta_{2}|\geq\epsilon}1\nonumber\\
&\leq & \mathcal{C} \frac{T^{a}}{\delta_{1}}T^{-\frac{2}{1+\vartheta_{kk}(\epsilon)}}\sum_{0<m\delta_{2}\leq T^{a}}1\nonumber\\
&\leq & \mathcal{C} T^{2a-\frac{2}{1+\vartheta_{kk}(\epsilon)}}(\delta_{1}\delta_{2})^{-1}\\
&= & \mathcal{C} T^{a-1} T^{a-\frac{1-\vartheta_{_{kk}}(\epsilon)}{1+\vartheta_{kk}(\epsilon)}}(\delta_{1}\delta_{2})^{-1}.
\end{eqnarray*}
Using the fact that (\ref{bkk}) and  $\limit{T} \delta_{i}(\ln T)^{1/\alpha}=\infty$ again, we have
$ S_{T,2}=o(T^{a-1})$ as $T\rightarrow\infty$.
}

The next lemma extends Lemma 2 of \cite{Pit2004} to the non-uniform sparse grid.
Let $\mathcal{R}(\delta)=\{t_{1}(T)<t_{2}(T)<....\}$ be a non-uniform grid such that
$$\delta_{max}:=\max_{t_{k}(T)\in [0,T]}(t_{k}(T)-t_{k-1}(T))\leq \delta_{0}\ \ \mbox{and}\ \ \ \delta_{min}(\ln T)^{1/\alpha}:=\min_{t_{k}(T)\in [0,T]}(t_{k}(T)-t_{k-1}(T))(\ln T)^{1/\alpha}\rightarrow\infty$$
as $T\rightarrow\infty$. 

\BL\label{lem:A5} For $S=T^a, a\in (0,1)$ we have for any $k\le p$ 
\begin{eqnarray*}
\mathbb{P}\left\{\max_{t\in \mathcal{R}(\delta)\cap[0,S]}\eta_k(t)> u_T\right\}=\sharp(\mathcal{R}(\delta)\cap[0,S])\overline{\Phi}(u_T)(1+o(1))
\end{eqnarray*}
as $u_T\to \IF $, where $\sharp(A)$ denotes the number of the elements  of the set $A$.
\EL

\prooflem{lem:A5} By  Bonferroni inequality for all $T$ large (set  $\Theta_{T}:=\sharp(\mathcal{R}(\delta)\cap[0,S])$  and $u:=u_T$)
\begin{eqnarray*}
\sum_{l=1}^{\Theta_{T}}\mathbb{P}\left\{\eta_{k}(t_{l}(T))> u\right\}&\geq& \mathbb{P}\left\{\max_{t\in \mathcal{R}(\delta)\cap[0,S]}\eta_{k}(t)> u\right\}\\
&\geq& \sum_{l=1}^{\Theta_{T}}\mathbb{P}\left\{\eta_{k}(t_{l}(T))> u\right\}-\sum_{1\leq m<l\leq \Theta_{T}}\mathbb{P}\left\{\eta_{k}(t_{m}(T))> u, \eta_{k}(t_{l}(T))> u\right\}\\
&=:& P_{1}+P_{2}.
\end{eqnarray*}
By the stationarity of $\eta_k$
$$P_{1}=\Theta_{T}\overline{\Phi}(u),$$
whereas for the second term we have for all $\ve>0$ sufficiently small
\begin{eqnarray*}
P_{2}&=&
\sum_{1\leq m<l\leq \Theta_{T}\atop t_{l}(T)-t_{m}(T)\leq \epsilon}\mathbb{P}\left\{\eta_{k}(t_{m}(T))> u, \eta_{k}(t_{l}(T))> u\right\}+\sum_{1\leq m<l\leq \Theta_{T}\atop t_{l}(T)-t_{m}(T)>\epsilon}\mathbb{P}\left\{\eta_{k}(t_{m}(T))> u, \eta_{k}(t_{l}(T))> u\right\}\\
&=:& P_{21}+P_{22}.
\end{eqnarray*}
Similarly as in the  calculations of  $R_{T,1}$ setting $r=r(t_{l}(T)-t_{m}(T))$ we have 
\begin{eqnarray*}
P_{21}&=&\sum_{1\leq m<l\leq \Theta_{T}\atop t_{l}(T)-t_{m}(T)\leq \epsilon}\mathbb{P}\left\{\eta_{k}(0)> u, \eta_{k}(t_{l}(T)-t_{m}(T))> u\right\}\\
&\leq&\sum_{1\leq m<l\leq \Theta_{T}\atop t_{l}(T)-t_{m}(T)\leq \epsilon}\overline{\Phi}(u)\overline{\Phi}\left(u\frac{\sqrt{1-r}}{1+r}\right).
\end{eqnarray*}
Since by (\ref{eq.A14})
$$\frac{1-r}{1+r}=\frac{1-r_{kk}(t_{l}(T)-t_{m}(T))}{1+r_{kk}(t_{l}(T)-t_{m}(T))}\geq \frac{1}{4}(t_{l}(T)-t_{m}(T))^{\alpha}\geq \frac{1}{4}\delta_{min}^{\alpha}$$
and the fact that  $u=u_T=(2\ln T)^{1/2}(1+o(1))$ and $\delta_{min}(\ln T)^{1/\alpha}\rightarrow\infty$ we get
\begin{eqnarray*}
P_{21}&\leq &\Theta_{T} \overline{\Phi}(u) \frac{\epsilon}{\delta_{min}} \overline{\Phi}\left(\frac{1}{2}u\delta_{min}^{\alpha/2}\right)\\
&\leq&\mathcal{C} \Theta_{T} \overline{\Phi}(u) \frac{\epsilon}{\delta_{min}} \frac{1}{u\delta_{min}^{\alpha/2}}\exp(-\frac{1}{8}u^{2}\delta_{min}^{\alpha})\\
&\leq&\mathcal{C}\Theta_{T} \overline{\Phi}(u)\frac{1}{u\delta_{min}^{\alpha/2}}\\
&=&\Theta_{T} \overline{\Phi}(u) o(1).
\end{eqnarray*}
Recalling the bound derived for $S_{T,2}$, by stationarity and Berman's inequality
\begin{eqnarray*}
P_{22}&\leq&\sum_{1\leq m<l\leq \Theta_{T}\atop t_{l}(T)-t_{m}(T)> \epsilon}\bigg[\overline{\Phi}^{2}(u)
+\mathcal{C}\exp\left(\frac{u^{2}}{1+\vartheta_{kk}(t_{l}(T)-t_{m}(T))}\right)\bigg]\\
&\leq& \mathcal{C}\Theta_{T}^{2}\bigg[\overline{\Phi}^{2}(u)
+\mathcal{C}\exp\left(\frac{u^{2}}{1+\vartheta_{kk}(\epsilon)}\right)\bigg]\\
\end{eqnarray*}
Noting that $\Theta_{T}\leq S/\delta_{min}=T^{a}/\delta_{min}$ and by repeating the calculations for $S_{T,2}$ we obtain further
$P_{22}=\Theta_{T} \overline{\Phi}(u) o(1)$ as $T\rightarrow\infty$, which completes the proof.
\QED

\BL\label{lem:A2} \jE{If}
$\mathfrak{R}(\delta_{1})$ and $\mathfrak{R}(\delta_{2})$ are sparse  and Pickands girds, respectively, then \Tan{for $k\leq p$} as $T\rightarrow\infty$
\begin{eqnarray*}
\mathbb{P}\left\{\max_{t\in \mathcal{R}(\delta_{1})\cap[0,S]}\eta_{k}(t)> u(y_{k1},z_{k}), \max_{t\in \mathcal{R}(\delta_{2})\cap[0,S]}\eta_{k}(t)> u(y_{k2},z_{k})\right\}=o(\LE{T^{a-1}}).
\end{eqnarray*}
\EL

\prooflem{lem:A2} Since $\mathfrak{R}(\delta_{1})$ and $\mathfrak{R}(\delta_{2})$ are sparse and Pickands girds, respectively, we have
$$\lim_{T\rightarrow\infty}\delta_{1}(T)/\delta_{2}(T)=\infty.$$
\Quan{Consequently}, the proof is similar to that of the case that  $\lim_{T\rightarrow\infty}\delta_{1}(T)/\delta_{2}(T)=\infty$ of \nelem{lem:A1}, and
therefore we omit further details. \QED

\aE{Let $X$ be a centered stationary Gaussian \Quan{process which satisfies condition \eqref{eq1.1} (as in the Introduction)}.
For the proof \Quan{of} Theorem 2.2  we shall determine the asymptotic \Quan{behaviours, as $u\to \IF$, of the} \jE{following probabilities}}
$$P_{S}(u,x)=\pk{\max_{t\in \mathfrak{R}(\delta_{1})\cap [0,S]}\aE{X(t)}
>u,\max_{t\in \mathfrak{R}(\delta_{2})\cap [0,S]}\aE{X(t)}>u+\frac{x}{u}}$$
and
$$P_{S}(u,x,y)=\pk{\max_{t\in \mathfrak{R}(\delta_{1})\cap [0,S]}\aE{X(t)}>u,\max_{t\in \mathfrak{R}(\delta_{2})\cap [0,S]}
\aE{X(t)}>u+\frac{x}{u},\max_{t\in [0,S]}\aE{X(t)}>u+\frac{y}{u}},$$
where $\mathfrak{R(\delta_{1})}=\mathfrak{R}(c(2\ln  T)^{-1/\alpha})$ and
$\mathfrak{R(\delta_{2})}=\mathfrak{R}(d(2\ln  T)^{-1/\alpha})$ with $c>d>0$.\\
For $\lambda\in ( c,\IF)$ along the lines of the proof
of Lemma D.1 \aE{in \cite{Pit96}} (see also the proof of Lemma 12.2.3 of \cite{leadbetter1983extremes}) 
$$P_{\lambda u^{-2/\alpha}}(u,x)\sim H_{c,d,\alpha}^{0,x}\overline{\Phi}(u)\ \
\mbox{and}\ \ P_{\lambda u^{-2/\alpha}}(u,x,y)\sim H_{c,d,\alpha}^{0,x,y}\overline{\Phi}(u)$$
as $u\rightarrow\infty$, where
%
$$\aE{H_{c,d,\alpha}^{0,x}}= \lim_{\lambda\to \IF}\frac{1}{\lambda}\int_{s\inr}e^{s}\pk{
\max_{\AH{k\inn:}kc\in[0,\lambda]}B^*_{\alpha/2}(kc)>s,
\max_{k\inn:kd\in [0,\lambda]}B^*_{\alpha/2}(kd)>s+x }\, ds$$
and
\begin{eqnarray*}
\aE{H_{c,d,\alpha}^{0,x,y}}&=& \lim_{\lambda\to \IF}\frac{1}{\lambda}\int_{s\inr}e^{s}\mathbb{P}\bigg\{
\max_{\AH{k\inn:}kc\in[0,\lambda]}B^*_{\alpha/2}(kc)>s,\\
&&\max_{k:kd\in [0,\lambda]}B^*_{\alpha/2}(kd)>s+x, \max_{\Quan{t}\in [0,\lambda]}
B^*_{\alpha/2}(\Quan{t})>s+y \bigg\}\, ds.
\end{eqnarray*}
\Quan{The next result can be shown \iH{along the same lines of} the proof of Theorem D.2 in  \cite{Pit96}.}

\BL\label{lem:A3} For any $x,y\inr$ we have
$$0<H_{c,d,\alpha}^{0,x}= \lim_{\lambda\to \IF}\frac{H_{c,d,\alpha}^{0,x}(\lambda)}{\lambda}<\infty\ \ \mbox{and}\ \
0<H_{c,d,\alpha}^{0,x,y}= \lim_{\lambda\to \IF}\frac{H_{c,d,\alpha}^{0,x,y}(\lambda)}{\lambda}<\infty.$$
Furthermore, \aE{for any $S>0$}
$$P_{S}(u,x)\sim S H_{c,d,\alpha}^{0,x}u^{2/\alpha}\overline{\Phi}(u)\ \ \mbox{and}\ \ P_{S}(u,x,y)\sim S H_{c,d,\alpha}^{0,x,y}u^{2/\alpha}\overline{\Phi}(u)$$
as $u\rightarrow\infty$.
\EL

\COM{Consequently, we have
\begin{eqnarray*}
&&P_{S}(u,x)=P\left(\max_{t\in \mathfrak{R}(\delta_{1})\cap [0,S]}X_{t}> u,\max_{t\in \mathfrak{R}(\delta_{2})\cap [0,S]}X_{t}> u+\frac{x}{u}\right)
\sim SH_{c,d,\alpha}^{0,x}u^{2/\alpha}\Psi(u).
\end{eqnarray*}
and
\begin{eqnarray*}
P_{S}(u,x,y)&=&P\left(\max_{t\in \mathfrak{R}(\delta_{1})\cap [0,S]}X_{t}>u,\max_{t\in \mathfrak{R}(\delta_{2})\cap [0,S]}X_{t}>u+\frac{x}{u},\max_{t\in \mathfrak{R}(\delta_{2})\cap [0,S]}X_{t}>u+\frac{y}{u}\right)\\
&\sim& SH_{c,d,\alpha}^{0,x,y}u^{2/\alpha}\Psi(u).
\end{eqnarray*}
}
%

\BL\label{lem:A4}
Let $\{Z_{T,ij},1 \le i \le p, 1\le j \le m\}, T>0$ be a random matrix. Suppose that the following convergence in distribution
$$\vk{Z}_{T,j}:=(Z_{T,1j} \ldot Z_{T,pj}) \todis (W_{1} \ldot W_p)=:\vk{W} , \quad T\to \IF$$
is valid for any index $j \le m$.  If further $Z_{T,ij} \le Z_{i1}$ \iH{holds} almost surely for any index $i\le p, 2 \le j\le m$, then we have the joint convergence in distribution
$$ ( \vk{Z}_{T,1} \ldot \vk{Z}_{T,k}) \todis (\vk{W} \ldot \vk{W}), \quad T\to \IF .$$
\EL
\prooflem{lem:A4} Assume for simplicity that $m=p=2$. By the assumptions,
Lemma 2.3 in \cite{JIH14} implies the convergence in distributions
$$ (Z_{T,11},Z_{T,12}) \todis (W_1,W_1), \quad (Z_{T,21},Z_{T,22}) \todis (W_2,W_2), \quad T\to \IF.$$
Hence we have the convergence in probability
$$ Z_{T,12}-Z_{T,11} \toprob 0, \quad Z_{T,22} -Z_{T,21}\toprob 0, \quad T\to \IF,$$
which then entails that 
$$ (Z_{T,11},Z_{T,21}, Z_{T,12},Z_{T,22}) \todis (W_1,W_2,W_1,W_2), \quad T\to \IF$$
establishing thus the proof. \QED

\textbf{Acknowledgments.}
\iH{E. Hashorva} kindly acknowledges partial support  by the Swiss National Science Foundation grant 200021-140633/1 and RARE -318984 (an FP7 Marie Curie IRSES Fellowship). Z. Tan acknowledges also support by the National Science Foundation of China (No. 11326175), RARE -318984 and Natural Science Foundation of Zhejiang Province of China (No. LQ14A010012).

\bibliographystyle{plain}

 \bibliography{difgridABCDEFG}
\end{document}